\documentclass[10pt]{article}
\usepackage[cp1252]{inputenc}
\usepackage[T1]{fontenc}
\usepackage{lmodern}
\usepackage[a4paper]{geometry}
\usepackage[english]{babel}
\usepackage{pifont}
\usepackage{extsizes}
\usepackage{xcolor}
\usepackage{fancybox}
\usepackage{graphicx}
\usepackage{sectsty}
\usepackage{enumerate}
\usepackage{multicol}
\usepackage[french]{varioref}
\usepackage[fleqn]{amsmath}
\usepackage{mathtools}
\usepackage{amssymb}
\usepackage{amsmath}
\usepackage{mathbbol}
\usepackage{mathrsfs}
\usepackage{bm}
\usepackage{amsthm}
\usepackage{blkarray}
\usepackage{subfigure}
\usepackage[babel=true]{csquotes}
\theoremstyle{break}
\newtheorem{Th}{Theorem}
\newtheorem{Prop}[Th]{Proposition}

\theoremstyle{definition}
\newtheorem{Rk}{Remark}

\newtheorem{Le}{Lemma}
\newtheorem{Cor}{Corollary}

\newcommand{\PP}{\mathbb{P}}
\newcommand{\EE}{\mathbb{E}}
\newcommand{\RR}{\mathbb{R}}

\newcommand{\diam}{\text{diam}}

\newcommand{\conv}[2][n]{\underset{#1\rightarrow #2}{\longrightarrow}}%écrit automatiquement la convergence
\newcommand{\EEE}[1]{\mathbb{E}\left[#1 \right]}%écrit l'espérance automatiquement avec les crochets
\newcommand{\PPP}[1]{\mathbb{P}\left(#1 \right)}%écrit la probabilité automatiquement avec la parenthèse
\newenvironment{prooft}[1]{\vskip 2mm\noindent {\bf Proof of #1.}}
                    {\hfill $\square$ \vskip 2mm \noindent}
\begin{document}
\pagenumbering{arabic}
\pagestyle{plain}
\author{\textbf{Pierre Calka}\footnote{Postal address: Université de Rouen, LMRS, avenue de l'Université, BP 12
76801 Saint-Etienne-du-Rouvray cedex, France. E-mail: pierre.calka@univ-rouen.fr} and \textbf{Nicolas Chenavier}\footnote{Postal address: Université de Rouen, LMRS, avenue de l'Université, BP 12
76801 Saint-Etienne-du-Rouvray cedex, France. E-mail: nicolas.chenavier@etu.univ-rouen.fr}}
\title{\textbf{Extreme values for characteristic radii of a Poisson-Voronoi tessellation}}
\maketitle

\begin{abstract}
A homogeneous Poisson-Voronoi tessellation of intensity $\gamma$ is observed in a convex body $W$. We associate to each cell of the tessellation two characteristic radii: the inradius, i.e. the radius of the largest ball centered at the nucleus and included in the cell, and the circumscribed radius, i.e. the radius of the smallest ball centered at the nucleus and containing the cell. We investigate the maximum and minimum of these two radii over all cells with nucleus in $W$. We prove that when $\gamma\rightarrow\infty$, these four quantities converge to Gumbel or Weibull distributions up to a rescaling. Moreover, the contribution of boundary cells is shown to be negligible. Such approach is motivated by the analysis of the global regularity of the tessellation. In particular, consequences of our study include the convergence to the simplex shape of the cell with smallest circumscribed radius and an upper-bound for the Hausdorff distance between $W$ and its so-called Poisson-Voronoi approximation.
\end{abstract}

\vspace{0.5cm}


\vspace{0.3cm}


%\makeatletter\tagsleft@true\makeatother
%\begin{align}
%\tag{H} a&=b
%\end{align}

\section{Introduction}
Let $\chi$ be a locally finite subset of $\RR^d$ endowed with its natural norm $|\cdot |$.  The Voronoi cell of nucleus $x\in \chi$ is the set \[C_{\chi}(x)=\{y\in\RR^d, |y-x|\leq |y-x'|, x\neq x'\in \chi\}.\] When $\chi=X_{\gamma}$ is a homogeneous Poisson point process of intensity $\gamma$, the family $\{C_{X_{\gamma}}(x), x\in X_{\gamma}\}$ is the so-called Poisson-Voronoi tessellation. Such model is extensively used in many domains such as cellular biology \cite{Po}, astrophysics \cite{RBFN}, telecommunications \cite{BB3} and ecology \cite{Ro}. For a complete account, we refer to the books \cite{OBSC}, \cite{SW}, \cite{M2} and the survey \cite{Cal5}. 

 To describe the mean behaviour of the tessellation, the notion of typical cell is introduced. The distribution of this random polytope can be defined as \[\EE[f(\mathscr{C}_{\gamma})] = \frac{1}{\gamma\lambda_d(B)}\EEE{\sum_{x\in X_{\gamma}\cap B}f(C_{X_{\gamma}}(x)-x)}\] where $f:\mathscr{K}^d\rightarrow\RR$ is any bounded measurable function on the set of convex bodies $\mathscr{K}^d$ (endowed with the Hausdorff topology), $\lambda_d$ is the $d$-dimensional  Lebesgue measure and $B$ is a Borel subset of $\RR^d$ with finite volume $\lambda_d(B)\in (0,\infty)$. Equivalently, $\mathscr{C}_{\gamma}$ is the Voronoi cell $C_{X_{\gamma}\cup\{0\}}(0)$ when we add the origin to the Poisson point process: this fact is a consequence of Slivnyak's Theorem, see e.g. Theorem 3.3.5  in \cite{SW}. The study of the typical cell in the literature includes mean values calculations \cite{M}, second order properties \cite{HM} and distributional estimates \cite{Cal6}, \cite{BL}, \cite{Mu}. A long standing conjecture due to D.G. Kendall about the asymptotic shape of large typical cell is proved in \cite{HRS2}.   

To the best of our knowledge, extremes of geometric characteristics of the cells, as opposed to their means, have not been studied in the literature up to now. In this paper, we are interested in the following problem: only a part of the tessellation is observed in a convex body $W$ (i.e. a convex compact set with non-empty interior) of volume $\lambda_d(W)=1$ where $\lambda_d$ denotes the Lebesgue measure in $\RR^d$. Let $f:\mathscr{K}^d\rightarrow\RR$ be a measurable function, e.g. the volume or the diameter of the cells. What is the limit behaviour of \[M_f(\gamma)=\max_{x\in X_{\gamma}\cap W}f(C_{X_{\gamma}}(x))\] when $\gamma$ goes to infinity? By scaling invariance of $X_{\gamma}$, it is the same as considering a tessellation with fixed intensity and observed in a window $W_{\rho}:=\rho W$ with $\rho\rightarrow\infty$. We give below some applications of such approach. 

First, the study of extremes describes the regularity of the tessellation. For instance, in finite element method, the quality of the approximation depends on some consistency measurements over the partition, see e.g. \cite{JGZ}. 

Another potential application field is statistics of point processes. The key idea would be to identify a point process from the extremes of its underlying Voronoi tessellation. A lot of inference methods have been developed for spatial point processes \cite{MW}. A comparison based on Voronoi extremes may or may not provide stronger results. At least, the regularity seems to discriminate to some extent to some point processes (see for instance a comparison between a determinantal point process and a Poisson point process in \cite{LH}). 

A third application is the so-called Poisson-Voronoi approximation i.e. a discretization of a convex body $W$ by the following union of Voronoi cells \[\mathscr{V}_{X_{\gamma}}(W)=\bigcup_{x\in X\cap W}C_{X_{\gamma}}(x).\] 
The first breakthrough is due to Heveling and Reitzner \cite{HR} and includes variance estimates of the volume of symmetric difference. However, the Hausdorff distance between the convex body and its approximation has not been studied yet. It is strongly connected to the maximum of the diameter of the cells which intersect the boundary of $\partial W$. We discuss this in section \ref{maxcirconsparagraph} and prove a rate of convergence of the approximation to the convex body with a suitable assumption on $W$. 
 
 Concretely, we are looking for two parameters $a_f(\gamma)$ and $b_f(\gamma)$ such that \[a_f(\gamma)M_f(\gamma)+b_f(\gamma)\underset{\gamma\rightarrow\infty}{\overset{\mathscr{D}}{\longrightarrow}}Y\] where $Y$ is a non degenerate random variable and $\overset{\mathscr{D}}{\longrightarrow}$ denotes the convergence in distribution. Up to a normalization, the extreme distributions of real random variables which are iid or with a mixing property are of three types: Fréchet, Gumbel or Weibull (see e.g. \cite{Lo} and \cite{L1}). More about extreme value theory can be found in the reference books by De Haan \& Ferreira \cite{HF} and by Resnick \cite{R}. Some extremes have been studied in stochastic geometry, for instance the maximum and minimum of inter-point distances of some point processes (see \cite{Pr}, \cite{MM} and \cite{JJ}) or the extremes of particular random fields \cite{LBB} but, to the best of our knowledge, nothing has been done for random tessellations. In our framework, the general theory cannot directly be applied for several reasons: unknown distribution of the characteristic for one fixed cell, dependency between cells and boundary effects. Moreover, the exceedances can be realized in clusters. For example, when the distance between the boundary of the cell and its nucleus is small, this is the same for one of its neighbors. Such clusters lead to the notion of extremal index, which was introduced by Leadbetter in \cite{L2}, and that we will study in a future work. 

 In this paper, we are interested in the characteristic radii i.e. inscribed and circumscribed radii of the Voronoi cell $C_{X_{\gamma}}(x)$ defined as \[r(C_{X_{\gamma}}(x))=\max\{r\geq 0, B(x,r)\subset C_{X_{\gamma}}(x)\} \text{ and } R(C_{X_{\gamma}}(x))=\min\{R\geq 0, B(x,R)\supset C_{X_{\gamma}}(x)\}\] where $B(x,r)$ is the ball of radius $r$ centered at $x$. Two reasons led us to the study of these quantities. First, the distribution tails of the inradius and circumscribed radius of the typical cell are easier to deal with \cite{Cal3} compared to other characteristics such as the volume or the number of hyperfaces. Secondly, knowing these two radii provides a better understanding of the cell shape since the boundary of $C_{X_{\gamma}}(x)$ is included in the annulus $B(x,R(C_{X_{\gamma}}(x)))-B(x,r(C_{X_{\gamma}}(x)))$. We consider the extremes 
\begin{equation}\label{def}\left.\begin{split} & r_{\max}(\gamma) = \max_{x\in X_{\gamma}\cap W}r(C_{X_{\gamma}}(x)),\hspace{0.5cm}  r_{\min}(\gamma) = \min_{x\in X_{\gamma}\cap W}r(C_{X_{\gamma}}(x))\\
& R_{\max}(\gamma) = \max_{x\in X_{\gamma}\cap W}R(C_{X_{\gamma}}(x)),\hspace{0.5cm} R_{\min}(\gamma) = \min_{x\in X_{\gamma}\cap W}R(C_{X_{\gamma}}(x)).\end{split}\right.  \end{equation}  In the following theorem, we derive the convergence in distribution of these quantities over cells with nucleus in $W$.

\begin{Th}
\label{Extr}
Let $X_{\gamma}$ be a Poisson point process of intensity $\gamma$ and $W$ a convex body of volume 1 in $\RR^d$. Then
\begin{subequations}
  \begin{equation}\label{maxins}
  \PP\left(2^d\kappa_d\gamma r_{\max}(\gamma)^d-\log(\gamma)\leq t \right) \underset{\gamma\rightarrow\infty}{\longrightarrow} e^{-e^{-t}}, \hspace{0.25cm} t\in\RR,\end{equation}
    \begin{equation}\label{minins}
  \PP\left(2^{d-1}\kappa_d\gamma^2r_{\min}(\gamma)^d\geq t \right) \underset{\gamma\rightarrow\infty}{\longrightarrow} e^{-t}, \hspace{0.25cm} t\geq 0,\end{equation}
  \begin{equation}\label{maxcircons}
  \PP\left(\kappa_d\gamma R_{\max}(\gamma)^d-\log\left(\alpha_1\gamma(\log\gamma)^{d-1}\right) \leq t\right) \underset{\gamma\rightarrow\infty}{\longrightarrow} e^{-e^{-t}}, \hspace{0.25cm} t\in\RR, \end{equation}
  \begin{equation}\label{mincircons}
  \PP(\alpha_2\kappa_d \gamma^{(d+2)/(d+1)}R_{\min}(\gamma)^d\geq t) \underset{\gamma\rightarrow\infty}{\longrightarrow} e^{-t^{d+1}}, \hspace{0.25cm} t\geq 0,\end{equation}
  \end{subequations}
  where $\alpha_1$ and $\alpha_2$ are given in \eqref{defalpha1} and \eqref{defalpha2} and $\kappa_d=\lambda_d(B(0,1))$. 
\end{Th}
The limit distributions are of type II and III and do not depend on the shape of $W$. One can note that the ratios  $r_{\max}(\gamma)/r_{\min}(\gamma)$ and $R_{\max}(\gamma)/R_{\min}(\gamma)$ are of respective orders  $(\gamma\log\gamma)^{1/d}$ and $(\gamma^{1/(d+1)}\log\gamma)^{1/d}$. This quantifies to some extent the irregularity of the Poisson-Voronoi tessellation. Moreover, the ratio $r_{\max}(\gamma)/R_{\max}(\gamma)$ is bounded. It suggests that  large cells tend to be spherical around the nucleus. This fact seems to confirm the D.G. Kendall's conjecture. 

As it is written, Theorem 1 is not applicable for concrete data. Indeed, in practice, the only cells which can be measured are included in the window. The following proposition addresses this problem. 

\begin{Prop}
\label{Boundarycharacteristic}
The extremes of characteristic radii over all cells included in $W$ or over all cells intersecting $\partial W$ have the same limit distributions as $r_{\max}(\gamma)$, $r_{\min}(\gamma)$, $R_{\max}(\gamma)$ and $R_{\min}(\gamma)$.
\end{Prop} 

\begin{center}
\begin{figure}[!h]
\label{densities}
\begin{center}
\begin{tabular}{cccc}
    \includegraphics[width=4.10cm,height=3.40cm]{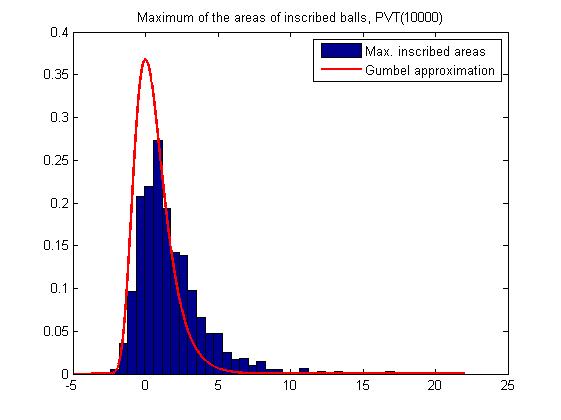} & \includegraphics[width=4.10cm,height=3.40cm]{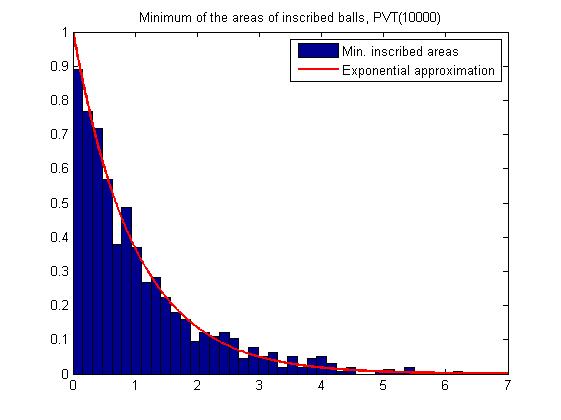}  & \includegraphics[width=4.10cm,height=3.40cm]{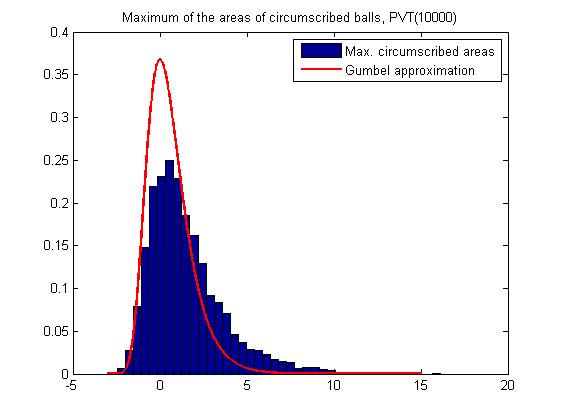}  &   \includegraphics[width=4.10cm,height=3.40cm]{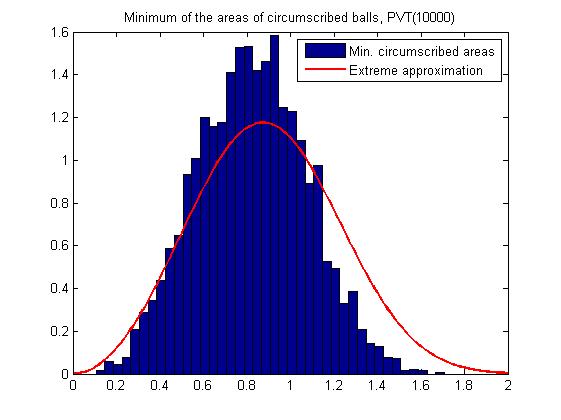}\\
     \includegraphics[width=4.10cm,height=3.40cm]{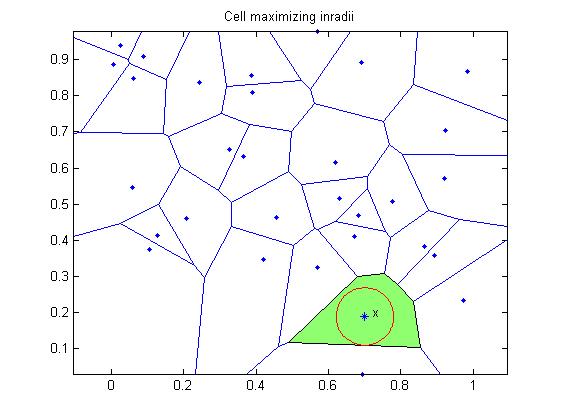} & \includegraphics[width=4.10cm,height=3.40cm]{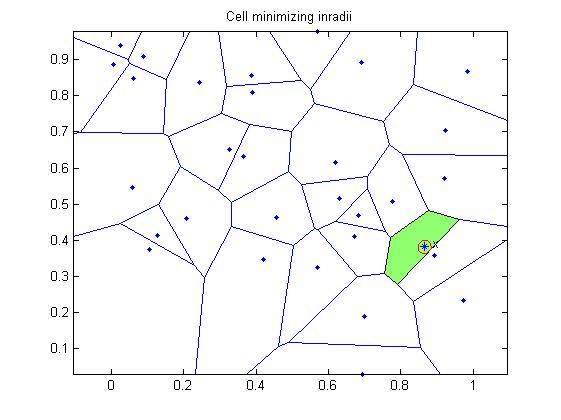} & \includegraphics[width=4.10cm,height=3.40cm]{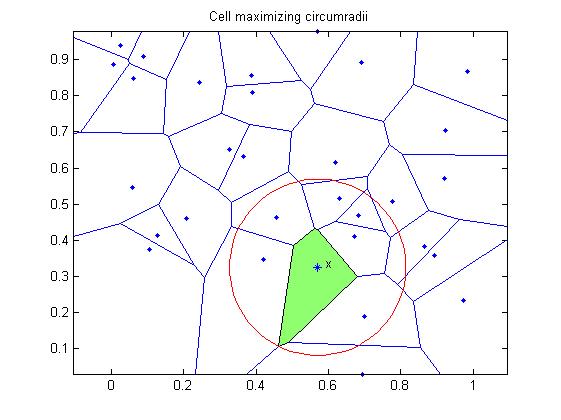}  & \includegraphics[width=4.10cm,height=3.40cm]{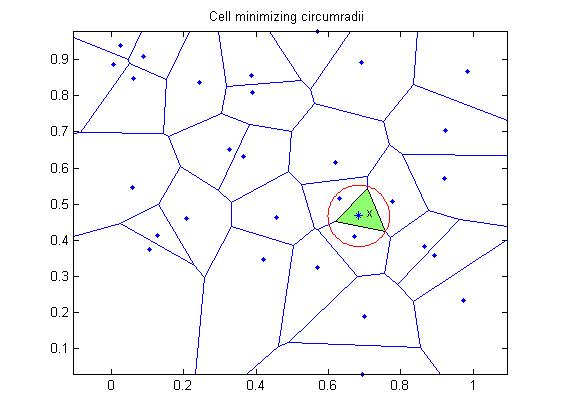}\\
     \textit{(a)} & \textit{(b)} & \textit{(c)} & \textit{(d)}
\end{tabular}
\end{center}
\caption{Empirical densities of the extremes based on 3500 simulations of PVT in $2\text{D}$ with $\gamma=10000$, for the cells included in $W=[0,1]^2$, on ${\text{Matlab}}^\copyright$. \textit{(a)} Cell maximizing the inradius. \textit{(b)} Cell minimizing the inradius. \textit{(c)} Cell maximizing the circumradius. \textit{(d)} Cell minimizing the circumradius.} 
\end{figure}
\end{center}

The convergences are illustrated in Figure \ref{densities} for the cells which are included in $W=[0,1]^2$. For sake of simplicity, the Poisson point process has been realized only in $W$. Because of Proposition \ref{Boundarycharacteristic} and related arguments, this does not affect the distribution over cells included in $W$. Simulations suggest that the rates of convergence are not the same for all these quantities. Indeed, in a future work, we will show that the rate is of the order of $\gamma^{-1}$, $\gamma^{-1/4}$ and $\gamma^{-1/6}$ for $r_{\min}(\gamma)$, $r_{\max}(\gamma)$ and $R_{\min}(\gamma)$ respectively.

All results of Theorem \ref{Extr} use geometric interpretations. For the circumscribed radii $R_{\max}(\gamma)$ and $R_{\min}(\gamma)$, we write the distributions as covering probabilities of spheres. The inscribed radii can be interpreted as interpoint distances. A study of the extremes of these distances has been done in several works such as \cite{JJ} and \cite{He}. For sake of completeness, we have rewritten these results in our setting in particular because the boundary effects are highly non trivial. Convergences \eqref{maxins} and \eqref{mincircons} could be obtained by considering underlying random fields and using methods inherited  for \cite{AGG} and \cite{S}. However, this approach does not provide \eqref{minins} and \eqref{maxcircons}. We will develop this idea in a future work and deduce some rates of convergence.

The paper is organized as follows. In section 2, we provide some preliminary result which shows that the boundary cells are negligible and implies Proposition \ref{Boundarycharacteristic}. In sections 3, 4 and 5, proofs of \eqref{mincircons}, \eqref{maxins}, \eqref{maxcircons} and \eqref{minins} are respectively given. Section 3 requires a technical lemma about deterministic covering of the sphere by caps which is proved in appendix. Section 4 contains an application of \eqref{maxcircons} to the Hausdorff distance between $W$ and its Poisson-Voronoi approximation. In section 5, we get a specific treatment of boundary effects which is more precise than in section 2.

In the rest of the paper, $c$ denotes a generic constant which does not depend on $\gamma$ but may depend on other quantities. The term $u_{\gamma}$ denotes a generic function of $t$, dependending on $\gamma$, which is specified at the beginning of sections 3,4 and 5. 

\section{Preliminaries on boundary effects}
\label{sectionboundary}
In this section, we show that the asymptotic behaviour of an extreme is in general not affected by  boundary cells. We  apply that result directly to the extremes of characteristic radii in order to show that Theorem \ref{Extr} implies Proposition \ref{Boundarycharacteristic}.

Let $f:\mathscr{K}_d\rightarrow\RR$ be a $k$-homogeneous measurable function, $0\leq k\leq d$ (i.e. $f(\lambda C)=\lambda^kf(C)$ for all $\lambda\in\RR_+$ and $C\in\mathscr{K}_d$). We consider for any $l\in\RR$
\begin{subequations}
  \begin{equation*}
 M_f^b(\gamma,l) = \max_{x\in X_{\gamma}, C_{X_{\gamma}}(x)\cap W_{1+l}\neq\varnothing}f(C_{X_{\gamma}}(x)),
  \end{equation*}
     \begin{equation*}
 M_f(\gamma,l) = \max_{x\in X_{\gamma}\cap W_{1+l}}f(C_{X_{\gamma}}(x)),
  \end{equation*}
   \begin{equation*}
 M_f^i(\gamma,l) = \max_{x\in X_{\gamma}, C_{X_{\gamma}}(x)\subset W_{1+l}}f(C_{X_{\gamma}}(x)),
  \end{equation*} 
  \end{subequations} where $W_{1+l}= (1+l)W$. When $l=0$, these maxima are simply denoted by $M_f^b(\gamma)$, $M_f(\gamma)$ and $M_f^i(\gamma)$. We define, for all $\epsilon >0$, a function $l_{\gamma}$ as  \begin{equation}\label{deflgamma}l_{\gamma}=\gamma^{-(1-\epsilon)/d}.\end{equation} Under suitable conditions, the following proposition shows that $M^b_f(\gamma)$, $M_f(\gamma)$ and $M_f^i(\gamma)$ satisfy the same convergence in distribution.

\begin{Prop}
\label{boundaryLemma}
Let $Y$ be a random variable and $a_{\gamma}$, $b_{\gamma}$ two functions such that \begin{equation}\label{hypothesisab}\frac{a_{\gamma}}{a_{\gamma_{\pm}}}\conv[\gamma]{\infty}1,\hspace{0.25cm} l_{\gamma}b_{\gamma}\conv[\gamma]{\infty}0 \text{ and } \frac{b_{\gamma}a_{\gamma_{\pm}}-a_{\gamma}b_{\gamma_{\pm}}}{a_{\gamma}}\conv[\gamma]{\infty}0\end{equation} with $\gamma_{+}=(1+l_{\gamma})^k\gamma$ and $\gamma_{-}=(1-l_{\gamma})^k\gamma$ for a certain $\epsilon$.  Then
\[a_{\gamma}M_f^b(\gamma)+b_{\gamma}\overset{\mathscr{D}}{\underset{\gamma\rightarrow\infty}{\longrightarrow}}Y  \Longleftrightarrow a_{\gamma}M_f(\gamma)+b_{\gamma}\overset{\mathscr{D}}{\underset{\gamma\rightarrow\infty}{\longrightarrow}}Y \Longleftrightarrow a_{\gamma}M_f^i(\gamma)+b_{\gamma}\overset{\mathscr{D}}{\underset{\gamma\rightarrow\infty}{\longrightarrow}}Y.\]

\end{Prop}

Before proving Proposition \ref{boundaryLemma}, we need an intermediary result due to Heinrich, Schmidt and Schmidt (Lemma 4.1 of \cite{H1}) which shows that, with high probability, the cells which intersect $\partial W$ have nucleus close to $\partial W$. Actually, they showed it for any stationary tessellation of intensity 1 which is observed in a window $\rho W$ with $\rho\rightarrow\infty$. For sake of completeness, we rewrite their result in a more explicit version for a Poisson-Voronoi tessellation. 
\begin{Le}(Heinrich, Schmidt and Schmidt)
\label{HSS}
Let us denote by $A_{\gamma}$ and $B_{\gamma}$ the events \[A_{\gamma}= \left\{\bigcap_{x\in X_{\gamma}}\{C_{X_{\gamma}}(x)\cap W=\varnothing\}\cup\{x\in W_{1+l_{\gamma}}\}\right\}\text{ and } B_{\gamma}= \left\{\bigcap_{x\in X_{\gamma}}\left\{C_{X_{\gamma}}(x)\subset W\right\}\cup\{x\not\in W_{1-l_{\gamma}}\}\right\} \] where $l_{\gamma}$ is given in \eqref{deflgamma}. Then $\PP(A_{\gamma})$ and $\PP(B_{\gamma})$ converge to 1 as $\gamma$ goes to infinity.\end{Le}
\begin{prooft}{Lemma \ref{HSS}}
In \cite{H1}, it is shown that \begin{equation}\label{HSSorigin}\PPP{\left\{\bigcap_{x\in X_{1}}\{C_{X_{1}}(x)\cap W_{\rho}=\varnothing\}\cup\{x\in W_{\rho+q(\rho)}\}\right\}\cap \left\{\bigcap_{x\in X_{1}}\left\{C_{X_{1}}(x)\subset W_{\rho}\right\}\cup\{x\not\in W_{\rho-q(\rho)}\}\right\}}\conv[\gamma]{\infty}1 \end{equation}  where $q(\rho)$ is the solution of the functional equation \[\rho^d=H(q^d(\rho)).\] The function $H:\RR_+\rightarrow\RR_+$ is convex, strictly increasing 
on its support $(x_o, \infty)$ (for some $x_0 \geq 0)$, such that $H(x)/x$ is non-decreasing for $x > 0$,
$\lim H(x)/x=\infty$ and $\EE[H (D^d(\mathscr{C}_1))]<\infty$ where $D(\mathscr{C}_1)$ is the diameter of the typical cell. 

In the case of a Poisson-Voronoi tessellation, $q(\rho)$ can be made explicit. Indeed, we can show that all moments of  $D(\mathscr{C}_1)$ exist since $D(\mathscr{C}_1)\leq 2R(\mathscr{C}_1)$ and $R(\mathscr{C}_1)$ is shown to have an exponentially decreasing tail in any dimension by an argument similar to Lemma 1 of \cite{Z}. Consequently, for a fixed $\epsilon\in (0,1)$, the functions $H$ and $q$ can be chosen as $H(x)=x^{1/\epsilon}$  and $q(\rho)=\rho^{\epsilon}$. 
Using the scaling property of Poisson point process, \[X_{1}\cap W_{\rho} \overset{\mathscr{D}}{=}\gamma^{1/d}(X_{\gamma}\cap W)\text{ and } X_{1}\cap W_{\rho\pm q(\rho)}\overset{\mathscr{D}}{=}\gamma^{1/d}(X_{\gamma}\cap W_{1\pm l_{\gamma}}) \] where $\gamma=\rho^d$ and $l_{\gamma}=\gamma^{-(1-\epsilon)/d}$. We deduce Lemma \ref{HSS} from \eqref{HSSorigin}. 
\end{prooft}

\begin{prooft}{Proposition \ref{boundaryLemma}}
\textit{First equivalence:} Let us assume that $a_{\gamma}M_f^b(\gamma)+b_{\gamma}$ converges in distribution to Y. On the event $A_{\gamma}$, $\forall x\in X_{\gamma}$, $C_{X_{\gamma}}(x)\cap W\neq\varnothing\Longrightarrow x\in W_{1+l_{\gamma}}$. Hence  \begin{equation}\label{lowuppmaxcircons}M_f^b(\gamma) \leq M_f(\gamma,l_{\gamma}) \leq M_f^b(\gamma,l_{\gamma})
.\end{equation} Because of Lemma \ref{HSS}, it is enough to show the convergence in distribution of the random variables conditionally on $A_{\gamma}$. Thanks to the scaling property of Poisson point process and the $k$-homogeneity of $f$ 
\begin{equation}\label{equaldistribution}
M_f^b(\gamma,l_{\gamma}) \overset{\mathscr{D}}{=}(1+l_{\gamma})^kM_f^b(\gamma_+)
\end{equation} with $\gamma_+=(1+l_{\gamma})^k\gamma$. We deduce from \eqref{lowuppmaxcircons},\eqref{equaldistribution} and \eqref{hypothesisab} that $a_{\gamma}M_f(\gamma,l_{\gamma})+b_{\gamma}$ converges in distribution to $Y$. By the scaling property, we get
\begin{equation}
\label{maxcircons4}
a_{\gamma}M_f(\gamma)+b_{\gamma}\overset{\mathscr{D}}{\underset{\gamma\rightarrow\infty}{\longrightarrow}}Y
.\end{equation}
Conversely, if \eqref{maxcircons4} holds then, using the fact that \[M_f(\gamma)\leq M_f^b(\gamma)\leq M_f(\gamma,l_{\gamma})\overset{\mathscr{D}}{=}(1+l_{\gamma})M_f(\gamma_+) \] and proceeding along the same lines, we get $a_{\gamma}M_f^b(\gamma)+b_{\gamma}\overset{\mathscr{D}}{\underset{\gamma\rightarrow\infty}{\longrightarrow}}Y$.

\textit{Second equivalence:} On the event $B_{\gamma}$, $\forall x\in X_{\gamma}$, $x\in W_{1-l_{\gamma}}\Longrightarrow C_{X_{\gamma}}(x)\subset W$. We prove the second equivalence as previously noting that, conditionally on $B_{\gamma}$ \begin{equation}\label{lowuppmaxcircons2}
M_f^i(\gamma,-l_{\gamma})\leq M_f(\gamma,-l_{\gamma})\leq M_f^i(\gamma)\leq M_f(\gamma).\end{equation} 
\end{prooft}

\section{Proof of \eqref{mincircons} and \eqref{maxins} }

\begin{prooft}{\eqref{mincircons}} Let $t\geq 0$ be fixed. We denote by $u_{\gamma}$ the following function: \begin{equation}\label{defugamma}u_{\gamma}=u_{\gamma}(t) = \left(\alpha_2^{-1}\kappa_d^{-1}\gamma^{-(d+2)/(d+1)}t \right)^{1/d}\end{equation} where $\alpha_2$ is given by \eqref{defalpha2}. Our aim is to prove that $\PP(R_{\min}(\gamma)\geq u_{\gamma})$ converges to $e^{-t^{d+1}}$ where $R_{\min}(\gamma)$ has been defined in \eqref{def}. The main idea is to deduce the asymptotic behaviour of $R_{\min}(\gamma)$ from the study of finite dimensional distributions $(R(C_{X_{\gamma}\cup\{\mathbf{x}_K\}}(x_1)),\ldots, R(C_{X_{\gamma}\cup\{\mathbf{x}_K\}}(x_K))$ for all $\{\mathbf{x}_K\} = \{x_1\ldots, x_K\}$ and $K\geq 1$. To do this, we write a new adapted version of a lemma due to Henze (see Lemma p. 345 in \cite{He}) in a context of  Poisson point process.

\begin{Le}
\label{rewriteHenze}
Let $f:\mathscr{K}_d\rightarrow\RR$, $F:\mathscr{K}_d\rightarrow\RR$ be two measurable functions and $A$ a Borel subset of $\RR$. Let us assume that for any $K\geq 1$, \begin{equation}\label{condHenze}\gamma^K\int_{W^K}\PPP{\forall i\leq K, f(C_{X_{\gamma}\cup\{\mathbf{x}_K\}}(x_i))<u_{\gamma}, F(C_{X_{\gamma}\cup\{\mathbf{x}_K\}}(x_i))\in A}d\mathbf{x}_K \conv[\gamma]{\infty}\lambda^K \end{equation} where  $d\mathbf{x}_K=dx_1\cdots dx_K$. Then
\[\PPP{\min_{x\in X_{\gamma}\cap W, F(C_{X_{\gamma}}(x))\in A}f(C_{X_{\gamma}}(x))\geq u_{\gamma}} \conv[\gamma]{\infty}e^{-\lambda}.\] 
\end{Le}

\begin{prooft}{Lemma \ref{rewriteHenze}}
Let $K$ be a fixed integer. The proof is close to the proof of Henze's Lemma and uses Bonferroni inequalities: one can show that if $A_{x,X_{\gamma}}$ is an $X_{\gamma}$-measurable event for all $x\in X_{\gamma}\cap W$, then  \begin{multline}\label{Bonferroni}\sum_{k=0}^{2K}\frac{(-1)^{k+1}}{k!}\EEE{\sum_{(x_1,\ldots, x_k)_{\neq}\in X_{\gamma}\cap W}\mathbb{1}_{A_{x_1,X_{\gamma}}}\ldots \mathbb{1}_{A_{x_k,X_{\gamma}}}} \leq \PPP{\bigcup_{x\in X_{\gamma}\cap W}A_{x,X_{\gamma}}}\\
\leq \sum_{k=0}^{2K+1}\frac{(-1)^{k+1}}{k!}\EEE{\sum_{(x_1,\ldots, x_k)_{\neq}\in X_{\gamma}\cap W}\mathbb{1}_{A_{x_1,X_{\gamma}}}\ldots \mathbb{1}_{A_{x_k,X_{\gamma}}}} .\end{multline} where $(x_1,\ldots, x_k)_{\neq}$ means that $(x_1,\ldots, x_k)$ is a $k$-tuple of distinct points. Applying \eqref{Bonferroni} to \[A_{x,X_{\gamma}}=\{f(C_{X_{\gamma}}(x))<u_{\gamma}\}\cap \{F(C_{X_{\gamma}}(x))\in A\},\] from Slivnyak's formula (see Corollary 3.2.3 in \cite{SW}), we obtain \begin{multline*}
\sum_{k=0}^{2K+1}\frac{(-1)^k}{k!}\gamma^k\int_{W^k}\PPP{\forall i\leq K, f(C_{X_{\gamma}\cup\{\mathbf{x}_K\}}(x_i))<u_{\gamma}, F(C_{X_{\gamma}\cup\{\mathbf{x}_K\}}(x_i))\in A}d\mathbf{x}_k\\
\leq \PPP{\min_{x\in X_{\gamma}\cap W, F(C_{X_{\gamma}}(x))\in A}f(C_{X_{\gamma}}(x))\geq u_{\gamma}}\\
\leq \sum_{k=0}^{2K}\frac{(-1)^k}{k!}\gamma^K\int_{W^k}\PPP{\forall i\leq K, f(C_{X_{\gamma}\cup\{\mathbf{x}_K\}}(x_i))<u_{\gamma}, F(C_{X_{\gamma}\cup\{\mathbf{x}_K\}}(x_i))\in A}d\mathbf{x}_k
.\end{multline*}
From \eqref{condHenze}, we obtain \begin{multline*}\sum_{k=0}^{2K+1}\frac{(-1)^k}{k!}\lambda^k\leq \liminf_{\gamma\rightarrow\infty}\PPP{\min_{x\in X_{\gamma}\cap W, F(C_{X_{\gamma}}(x))\in A}f(C_{X_{\gamma}}(x))\geq u_{\gamma}}\\
\leq \limsup_{\gamma\rightarrow\infty}\PPP{\min_{x\in X_{\gamma}\cap W, F(C_{X_{\gamma}}(x))\in A}f(C_{X_{\gamma}}(x))\geq u_{\gamma}}\leq \sum_{k=0}^{2K}\frac{(-1)^k}{k!}\lambda^k.\end{multline*} We conclude the proof by taking $K\rightarrow\infty$. 
\end{prooft}

We apply Lemma \ref{rewriteHenze} to  $f(C_{X_{\gamma}}(x))=R(C_{X_{\gamma}}(x))$. The function $F(C_{X_{\gamma}}(x))=F_{d-1}(C_{X_{\gamma}}(x))$ denotes the number of hyperfaces of the cell $C_{X_{\gamma}}(x)$.  In all the proof, the event considered is $A=\RR$. We notice that the choice of the function $F$ is of no importance here but will be essential in the proof of Propositions \ref{minsimplex} and \ref{minsimplexDequals2}. From Lemma \ref{rewriteHenze}, it is sufficient to study the limit behaviour of 
\begin{equation}\label{slivnyak}\gamma^K\int_{W^K}\PPP{\forall i\leq K, R(C_{X\cup\{\mathbf{x}_K\}}(x_i))<u_{\gamma}}d\mathbf{x}_K\end{equation} for all integer $K$. We divide the proof into two parts. 

\paragraph{Step 1 }
When $K=1$, using the stationarity of $X_{\gamma}$ and the fact that $\lambda_d(W)=1$, we show that the integral \eqref{slivnyak} is $\gamma\PP(R(C_{X_{\gamma}\cup\{0\}}(0))<u_{\gamma})$. As in \cite{Cal5} section 5.2.3, we can reinterpret the distribution function of $R(C_{X_{\gamma}\cup\{0\}}(0))$ as a covering probability to get \begin{equation}\label{eqcovering}\gamma\PP(R(C_{X_{\gamma}\cup\{0\}}(0))<u_{\gamma}) = \gamma \sum_{k=0}^{\infty}e^{-2^d\kappa_d\gamma u_{\gamma}^d}\frac{(2^d\kappa_d\gamma u_{\gamma}^d)^k}{k!}p_k\end{equation} where $p_k$ is the probability to cover the unit sphere with $k$ independent spherical caps such that their normalized radii are distributed as $d\nu(\theta) = d\pi\sin(\pi\theta)\cos^{d-1}(\pi\theta)\mathbb{1}_{[0,1/2]}(\theta)d\theta$. The equality comes from the fact that \begin{equation*}\begin{split} R(C_{X_{\gamma}\cup\{0\}}(0))<u_{\gamma} & \Longleftrightarrow  \text{ the family } \{\mathscr{A}_{y}(0), y\in X_{\gamma}\} \text{ covers } S(0,u_{\gamma})\\
& \Longleftrightarrow \text{ the family }\{\mathscr{A}_{y}(0), y\in X_{\gamma}\cap B(0,2u_{\gamma})\} \text{ covers } S(0,u_{\gamma})\end{split}\end{equation*} 
where \begin{equation}\label{defsphericalcap}\mathscr{A}_{y}(x)=S(x,u_{\gamma})\cap  H_y^+(x)\end{equation} and $H_y^+(x)$ is the half-space which contains $y$ and delimited by the bisecting hyperplane of $[x,y]$.

\begin{figure}[h!]
\begin{center}
\includegraphics[width=8.5cm,height=7cm]{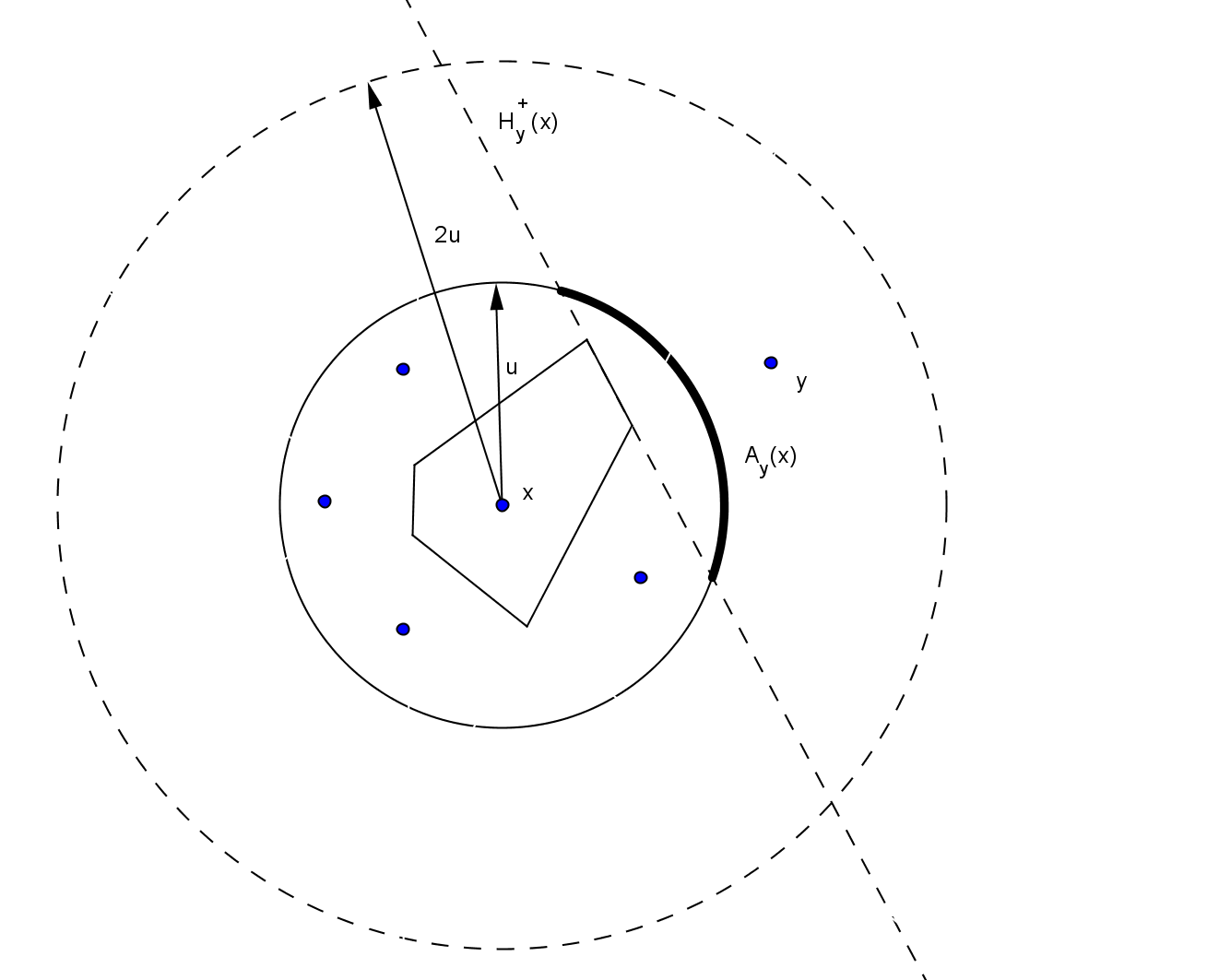}
  \caption{Interpretation of the circumscribed radius as a covering of sphere.}
\end{center}
\end{figure}
 We denote by \begin{equation}\label{defalpha2}
 \alpha_2 := \left(\frac{2^{d(d+1)}}{(d+1)!}p_{d+1}\right)^{1/(d+1)}>0.\end{equation} For example, when $d=2$, $\alpha_2=\left(\frac{5}{12}-\frac{4}{\pi^2} \right)^{1/3}$. 
 
 Since $p_k=0$ for all $k\leq d$, \eqref{eqcovering} gives \[\gamma\PP(R(C_{X_{\gamma}\cup\{0\}}(0))<u_{\gamma}) = \gamma \frac{(2^d\kappa_d\gamma u_{\gamma}^d)^{d+1}}{(d+1)!}e^{-2^d\kappa_d \gamma u_{\gamma}^d}p_{d+1} + \gamma\sum_{k=d+2}^{\infty}e^{-2^d\kappa_d\gamma u_{\gamma}^d}\frac{(2^d\kappa_d\gamma u_{\gamma}^d)^k}{k!}p_k.\] The first term converges to $t^{d+1}$ from \eqref{defugamma} and \eqref{defalpha2}. The second term is negligible since $\gamma (\gamma u_{\gamma}^d)^{d+2}=c\cdot \gamma^{-1/(d+1)}$ converges to 0 as $\gamma$ tends to infinity. This shows that \begin{equation}\label{cclstep1}\gamma\int_W\PP(R(C_{X_{\gamma}\cup\{x\}}(x))<u_{\gamma})dx \conv[\gamma]{\infty}t^{d+1}.\end{equation}

\paragraph{Step 2 }
When $K\geq 2$, we use the same interpretation as in step 1: for all $\mathbf{x}_K=(x_1,\ldots, x_K)\in W^K$, and $i\leq K$
\begin{equation*}\begin{split} R(C_{X_{\gamma}\cup\{\mathbf{x}_K\}}(x_i))<u_{\gamma} & \Longleftrightarrow \text{ the family } 
\{\mathscr{A}_y(x_i), y\in X_{\gamma}\cup\{\mathbf{x}_K\}-\{x_i\}\} \text{ covers } S(x_i,u_{\gamma})\\
& \Longleftrightarrow  \text{ the family } \{\mathscr{A}_y(x_i), y\in (X_{\gamma}\cup\{\mathbf{x}_K\}-\{x_i\})\cap B(x_i, 2u_{\gamma})\} \text{ covers } S(x_i,u_{\gamma}).
\end{split}\end{equation*} 
Hence, writing the previous event as \enquote{$S(x_i,u_{\gamma}) \text{ covered}$}, we have
 \begin{equation}\label{coveringproba}\PPP{\forall i\leq K, R(C_{X_{\gamma}\cup\{\mathbf{x}_K\}}(x_i))<u_{\gamma}} = \PPP{\bigcap_{i\leq K}\{ S(x_i,u_{\gamma}) \text{ covered}\}} 
 .\end{equation} 
We have now to consider the spherical caps induced by both the points $x_j, j\neq i$ and the points from $X_{\gamma}$.  For all $\mathbf{x}_K=(x_1,\ldots, x_K)\in W^K$, we denote by $n_l(\mathbf{x}_K)$ the number of connected components of $\bigcup_{i=1}^KB(x_i,2u_{\gamma})$ with exactly $l$ balls. Given $n_1,\ldots, n_K$  such that $\sum_{l=1}^{K}ln_l=K$, we define \begin{equation}\label{defW}W_{K}(n_1,\ldots, n_K) =\{\mathbf{x}_K\in W^K, n_l(\mathbf{x}_K)=n_l\text{ for all }l\leq K\}.\end{equation} Let us note that the subsets $W_{K}(n_1,\ldots, n_K)$, with $\sum_{l=1}^Kln_l=K$,  partition $W^K$. We then deal with two cases. 
\begin{enumerate}
\item If $B(x_i,2u_{\gamma})\cap B(x_j,2u_{\gamma})=\varnothing$ for all $i\neq j\leq K$ i.e. $\mathbf{x}_K\in W_{K}(K,\cdots, 0)$, the events considered in the right-hand side of \eqref{coveringproba} are independent.
\item If not, we are going to show that the contribution of such $\mathbf{x}_K$ in \eqref{slivnyak} is negligible. 
\end{enumerate}
%\begin{figure}[h!]
%\begin{center}
% \includegraphics[width=4.5cm,height=4cm]{Connectedcomponent.png}
 % \caption{Connected components of $\bigcup_{i=1}^KB(x_i,2u_{\gamma})$}
%\end{center}
%\end{figure}
More precisely, we write the integral \eqref{slivnyak} in the following way \begin{multline}\label{planintegral}
\gamma^K\int_{W^K}\PPP{\forall i\leq K, R(C_{X_{\gamma}\cup\{\mathbf{x}_K\}}(x_i))<u_{\gamma}}d\mathbf{x}_K = \gamma^K\int_{W_K(K,0,\ldots, 0)}\PPP{\forall i\leq K, R(C_{X_{\gamma}\cup\{\mathbf{x}_K\}}(x_i))<u_{\gamma}}d\mathbf{x}_K\\ + \gamma^K\int_{W^K-W_K(K,0,\ldots, 0)}\PPP{\forall i\leq K, R(C_{X_{\gamma}\cup\{\mathbf{x}_K\}}(x_i))<u_{\gamma}}d\mathbf{x}_K
.\end{multline}  

\subparagraph{Step 2.1}(Case of disjoint balls)
For all $\mathbf{x}_K=(x_1,\ldots, x_K)\in W_{K}(K,0,\ldots, 0)$, we obtain from \eqref{coveringproba} and \eqref{cclstep1}
\begin{equation}\label{disjointballs1}\begin{split}\gamma^K\PPP{\forall i\leq K, R(C_{X_{\gamma}\cup\{\mathbf{x}_K\}}(x_i))<u_{\gamma}} = \prod_{i=1}^K\gamma\PP(R(C_{X_{\gamma}\cup\{x_i\}}(x_i))<u_{\gamma}) \conv[\gamma]{\infty}(t^{d+1})^K \end{split}.\end{equation}  Moreover, $\lambda_{dK}(W_K(K,0,\ldots, 0)) \conv[\gamma]{\infty}1$. This shows that \begin{equation}\label{cclstep2.1}\gamma^K\int_{W_K(K,0,\ldots, 0)}\PPP{\forall i\leq K, R(C_{X_{\gamma}\cup\{\mathbf{x}_K\}}(x_i))<u_{\gamma}}d\mathbf{x}_K \conv[\gamma]{\infty}(t^{d+1})^K.\end{equation}

\subparagraph{Step 2.2}(Case of non disjoint balls)
In this step, we show that the second integral in the right-hand side of \eqref{planintegral} converges to 0. In particular, we study the limit behaviour of the integrand of \eqref{slivnyak} for all $\mathbf{x}_K=(x_1,\ldots, x_K)\in W_K(n_1,\ldots, n_K)$ with $(n_1,\ldots, n_K)\neq (K,0,\ldots, 0)$. The number of points of $X_{\gamma}\cap \bigcup_{i=1}^KB(x_i,2u_{\gamma})$ is Poisson distributed of mean $\gamma\lambda_d\left(\bigcup_{i=1}^KB(x_i,2u_{\gamma})\right)$. From \eqref{coveringproba}, we deduce that 
 \begin{multline}\label{singlecomponent}\gamma^K\PPP{\forall i\leq K, R(C_{X_{\gamma}\cup\{\mathbf{x}_K\}}(x_i))<u_{\gamma}}\\ = \gamma^K \sum_{k=0}^{\infty}\frac{\left(\gamma\lambda_d\left(\bigcup_{i=1}^kB(x_i,2u_{\gamma}) \right) \right)^k}{k!}e^{-\gamma\lambda_d\left(\bigcup_{i=1}^kB(x_i,2u_{\gamma}) \right)} \times p_k(x_1,\ldots, x_K).\end{multline}
The term $p_k(x_1,\ldots, x_K)$ denotes the probability to cover the spheres $S(x_i,u_{\gamma})$, $i=1\ldots K$, with the spherical caps $\{\mathscr{A}_{x_j}(x_i), i\neq j\leq K\}$ and $\{\mathscr{A}_{y_m}(x_i), m\leq k\}$, defined in \eqref{defsphericalcap}, where $y_1,\ldots, y_k$ are $k$ independent points which are uniformly distributed in $\bigcup_{i=1}^KB(x_i,2u_{\gamma})$. This probability satisfies the following property:
\begin{Le}
\label{Calotte}
Let $\mathbf{x}_K=(x_1,\ldots, x_K)\in W_K(n_1,\ldots, n_K)$ and \begin{equation}\label{defN}N=\sum_{l=1}^K(d+1)n_l.\end{equation} Then, for all $k<N$ 
 \begin{equation}\label{coverKspheres}p_k(x_1,\ldots, x_K)=0.\end{equation}  
\end{Le}

The proof of Lemma \ref{Calotte} is postponed to the appendix. From \eqref{singlecomponent}, \eqref{coverKspheres} and the trivial inequalities $0\leq p_k(x_1,\ldots, x_K)\leq 1$ and  $\lambda_d\left(\bigcup_{i=1}^kB(x_i,2u_{\gamma}) \right)\leq k2^d\kappa_du_{\gamma}^d$, we deduce that there exists a constant $c$, depending on $K$, such that
\[ \gamma^K\PPP{\forall i\leq K, R(C_{X_{\gamma}\cup\{\mathbf{x}_K\}}(x_i))<u_{\gamma}}  \leq  \gamma^K\sum_{k=N}^{\infty}\frac{(k2^d\kappa_d\gamma u_{\gamma}^d)^k}{k!} \underset{\gamma\rightarrow\infty}{\sim}c\cdot \gamma^K(\gamma u_{\gamma}^d)^{N}\] where $\phi(\gamma)\underset{\gamma\rightarrow\infty}{\sim}\psi(\gamma)$ means $\frac{\phi(\gamma)}{\psi(\gamma)}\conv[\gamma]{\infty}1$. Using \eqref{defugamma}, \eqref{defN} and the fact that $K=\sum_{l=1}^Kln_l$, we obtain for $\gamma$ large enough
\begin{equation}
\label{disjointballs2}
\gamma^K\PPP{\forall i\leq K, R(C_{X_{\gamma}\cup\{\mathbf{x}_K\}}(x_i))<u_{\gamma}}\leq c\cdot \prod_{l=2}^K\gamma^{(l-1)n_l}
.\end{equation}

Moreover, using the fact that $\lambda_{dK}(W_{K}(n_1,\ldots, n_K))\leq c\cdot\prod_{l=2}^K(u_{\gamma}^d)^{(l-1)n_l} = c\cdot \prod_{l=2}^K\gamma^{-\frac{(d+2)(l-1)}{d+1}n_l}$ and \eqref{disjointballs2}, we get
   \begin{multline}\label{majrest}\gamma^K\int_{W^{K}-W_K(K,0,\ldots, 0)}\PPP{\forall i\leq K, R(C_{X_{\gamma}\cup\{\mathbf{x}_K\}}(x_i))<u_{\gamma}}d\mathbf{x}_K\\
    = \sum \gamma^K\int_{W_{K}(n_1,\ldots, n_K)}\PPP{\forall i\leq K, R(C_{X_{\gamma}\cup\{\mathbf{x}_K\}}(x_i))<u_{\gamma}}d\mathbf{x}_K \leq c\cdot \sum\prod_{l=2}^K{\gamma}^{-\frac{l-1}{d+1}n_l} .\end{multline} The sum above runs over all the $K$-tuples $(n_1,\ldots,n_K)$ such that $\sum_{l=1}^Kln_l=K$ and $n_{1}\neq K$. Since $(n_1,\ldots, n_K)\neq (K,0,\ldots,0)$, there exists $l\geq 2$ such that $n_l\neq 0$. Consequently, we get from \eqref{majrest} \begin{equation}\label{cclstep2.2} \gamma^K\int_{W^{K}-W_K(K,0,\ldots, 0)}\PPP{\forall i\leq K, R(C_{X_{\gamma}\cup\{\mathbf{x}_K\}}(x_i))<u_{\gamma}}d\mathbf{x}_K =  O\left({\gamma}^{-1/(d+1)}\right)\end{equation} where $\phi(\gamma)=O(\psi(\gamma))$ means that $\frac{\phi(\gamma)}{\psi(\gamma)}$ is bounded. 
    
\paragraph{Conclusion}    
    
    From \eqref{cclstep2.1} and \eqref{cclstep2.2}, we deduce that for all $K\geq 1$ \[\gamma^K\int_{W^K}\PPP{\forall i\leq K, R(C_{X_{\gamma}\cup\{\mathbf{x}_K\}}(x_i))<u_{\gamma}}d\mathbf{x}_K\conv[\gamma]{\infty}(t^{d+1})^K.\] We then apply Lemma \ref{rewriteHenze}, with $A=\RR$, to conclude that \[\PPP{R_{\min}(\gamma)\geq u_{\gamma}}\conv[\gamma]{\infty}e^{-t^{d+1}}.\]
\end{prooft}

The cell which minimizes the circumscribed radius is asymptotically a simplex. To show it, we denote by \[R_{\min}'(\gamma) = \min_{x\in X_{\gamma}\cap W, F_{d-1}(C_{X_{\gamma}}(x))\geq d+2}R(C_{X_{\gamma}}(x))\] where  $F_{d-1}(C_{X_{\gamma}}(x))$ is the number of hyperfaces of $C_{X_{\gamma}}(x)$. The order of convergence of $R'_{\min}(\gamma)$ is greater than $u_{\gamma}$ according to the following proposition.
\begin{Prop}
\label{minsimplex}
Let $X_{\gamma}$ be a Poisson point process of intensity $\gamma$ and $W$ a convex body of volume 1. Then, for all $t\geq 0$, 
\[\PPP{\alpha_2\kappa_d\gamma^{(d+2)/(d+1)}R'^d_{\min}(\gamma)\geq t}\conv[\gamma]{\infty}1 .\] 
\end{Prop}
\begin{prooft}{Proposition \ref{minsimplex}}
We apply Lemma \ref{rewriteHenze} to $f(C_{X_{\gamma}}(x)) = R(C_{X_{\gamma}}(x))$ and $A=[d+2,\infty)$. We then study the finite dimensional distributions i.e. \begin{equation}
\label{condminsimplex}
\gamma^K\int_{W^K}\PPP{\forall i\leq K, R(C_{X_{\gamma}\cup\{\mathbf{x}_K\}}(x_i))<u_{\gamma}, F_{d-1}(C_{X_{\gamma}\cup\{\mathbf{x}_K\}}(x_i))\geq d+2}d\mathbf{x}_K
\end{equation}
for all $K\geq 1$. When $K=1$, the integrand of \eqref{condminsimplex} is \begin{multline*}\gamma\PPP{R(C_{X_{\gamma}\cup\{0\}}(0))<u_{\gamma}, F_{d-1}(C_{X_{\gamma}\cup\{0\}}(0))\geq d+2} \leq \gamma\PPP{R(C_{X_{\gamma}\cup\{0\}}(0))<u_{\gamma}, \#(X_{\gamma}\cap B(0,2u_{\gamma}))\geq d+2}\\ = \gamma\sum_{k=d+2}^{\infty}\frac{(2^d\kappa_d\gamma u_{\gamma}^d)^k}{k!}e^{-2^d\kappa_d\gamma u_{\gamma}^d}p_k\underset{\gamma\rightarrow\infty}{\sim}c\cdot\gamma^{-1/(d+1)}.\end{multline*}  We deduce that $\gamma\int_W\PPP{R(C_{X_{\gamma}\cup\{x\}}(x))<u_{\gamma}, F_{d-1}(C_{X_{\gamma}\cup\{x\}}(x))\geq d+2}$ converges to 0. More generally, for all $K\geq 1$, we get \begin{equation}
\label{condminsimplex1}
\gamma^K\int_{W_K(K,0,\ldots, 0)}\PPP{\forall i\leq K, R(C_{X_{\gamma}\cup\{\mathbf{x}_K\}}(x_i))<u_{\gamma}, F_{d-1}(C_{X_{\gamma}\cup\{\mathbf{x}_K\}}(x_i))\geq d+2}d\mathbf{x}_K\conv[\gamma]{\infty}0
.\end{equation}
Moreover, from \eqref{cclstep2.2} \begin{equation}
\label{condminsimplex2}
\gamma^K\int_{W^K-W_K(K,0,\ldots, 0)}\PPP{\forall i\leq K, R(C_{X_{\gamma}\cup\{\mathbf{x}_K\}}(x_i))<u_{\gamma}, F_{d-1}(C_{X_{\gamma}\cup\{\mathbf{x}_K\}}(x_i))\geq d+2}d\mathbf{x}_K\conv[\gamma]{\infty}0
.\end{equation}
From \eqref{condminsimplex1}, \eqref{condminsimplex2} and Lemma \ref{rewriteHenze} applied to $A=[d+2,\infty)$, we get \[\PPP{R'_{\min}(\gamma)\geq u_{\gamma}} \conv[\gamma]{\infty}1.\]
\end{prooft}

\begin{Cor}
\label{corollaryhausdorff}
Let $X_{\gamma}$ be a Poisson point process of intensity $\gamma$ and $W$ a convex body of volume 1. Then \[\PPP{\forall x\in X_{\gamma}, R(C_{X_{\gamma}}(x))=R_{\min}(\gamma)\Longrightarrow F_{d-1}(C_{X_{\gamma}}(x))=d+1}\conv[\gamma]{\infty}1.\]
\end{Cor} 

Proposition \ref{minsimplex} implies Corollary \ref{corollaryhausdorff} but does not provide the exact order of $R'_{\min}(\gamma)$. Nevertheless, when $d=2$, it can be made explicit. The key idea is contained in Lemma \ref{Calottesimplex} and cannot unfortunately be extended to higher dimensions. 

\begin{Prop}
\label{minsimplexDequals2}
Let $X_{\gamma}$ be a Poisson point process of intensity $\gamma$ and $W$ a convex body of volume 1 in $\RR^2$. Then, for all $t\geq 0$, 
\[\PPP{\alpha_2'\pi\gamma^{5/4}R'^2_{\min}(\gamma)\geq t}\conv[\gamma]{\infty}e^{-t^{4}} \] where $\alpha_2'$ is defined in \eqref{defalpha3}. 
\end{Prop}

\begin{prooft}{Proposition \ref{minsimplexDequals2}}
Let $t\geq 0$ be fixed and let us denote by \begin{equation}\label{defu'gamma}
u'_{\gamma}=u'_{\gamma}(t)=\left(\alpha_2'^{-1}\pi^{-1}\gamma^{-5/4}t \right)^{1/2}
\end{equation} where $\alpha_2'$ is specified in \eqref{defalpha3}. As in the proof of \eqref{mincircons}, we interpret the distribution function of $R'_{\min}(\gamma)$ as a covering probability of the circle. Let $\mu_k$ be the probability that $S(0,u'_{\gamma})$ is covered with the circular caps $\{\mathscr{A}_{y_m}(0), m\leq k\}$ where $y_1,\ldots, y_k$ are $k$ independent points which are uniformly distributed in $B(0,2u'_{\gamma})$ and such that $F_1(C_{\{0\}\cup\{\mathbf{y}_k\}}(0))\geq 4$ i.e. 
\begin{equation}
\label{muk}
\PPP{R(C_{X_{\gamma}\cup\{0\}}(0))<u'_{\gamma}, F_{1}(C_{X_{\gamma}\cup\{0\}}(0))\geq 4} = \sum_{k=4}^{\infty}\frac{1}{k!}(4\pi\gamma u'^2_{\gamma})^ke^{-4\pi\gamma u'^2_{\gamma}}\mu_k
.\end{equation}
The constant $\alpha_2'$ is defined as \begin{equation}\label{defalpha3}\alpha_2'=\left(\frac{32}{3}\mu_{4}\right)^{1/4}>0.\end{equation} We are going to apply Lemma \ref{rewriteHenze} to the event $A=[4,\infty)$ replacing $u_{\gamma}$ by $u'_{\gamma}$. To do it, we need to get the limit behaviour of \begin{equation}\label{aimsimplex}\gamma^K\int_{W^K}\PPP{\forall i\leq K, R(C_{X_{\gamma}\cup\{\mathbf{x}_K\}}(x_i))<u'_{\gamma}, F_{1}(C_{X_{\gamma}\cup\{\mathbf{x}_K\}}(x_i))\geq 4}d\mathbf{x}_K\end{equation} for all $K\geq 1$. 
 
  When $K=1$, from \eqref{muk} and \eqref{defu'gamma}, we deduce that $\gamma\int_W\PPP{ R(C_{X_{\gamma}\cup\{x\}}(x))<u'_{\gamma}, F_{d-1}(C_{X_{\gamma}\cup\{x\}}(x)))\geq 4}dx$ converges to $t^{4}$. More generally, for all $K\geq 1$, \begin{equation}\label{minsimplexstep2.1}
\gamma^K\int_{W_K(K,0,\ldots, 0)}\PPP{\forall i\leq K, R(C_{X_{\gamma}\cup\{\mathbf{x}_K\}}(x_i))<u'_{\gamma}, F_{1}(C_{X_{\gamma}\cup\{\mathbf{x}_K\}}(x_i))\geq 4}d\mathbf{x}_K\conv[\gamma]{\infty}t^{4K}
.\end{equation}

Otherwise, for all $\mathbf{x}_K\in W_K(n_1,\ldots, n_K)$ with $(n_1,\ldots, n_K)\neq (K,0,\ldots, 0)$, the integrand of \eqref{aimsimplex} is  \begin{multline}\label{jointdistributionsimplex}\gamma^K\PPP{\forall i\leq K, R(C_{X_{\gamma}\cup\{\mathbf{x}_K\}}(x_i))<u'_{\gamma}, F_{1}(C_{X_{\gamma}\cup\{\mathbf{x}_K\}}(x_i))\geq 4}\\ = \gamma^K \sum_{k=0}^{\infty}\frac{\left(\gamma\lambda_d\left(\bigcup_{i=1}^kB(x_i,2u'_{\gamma}) \right) \right)^k}{k!}e^{-\gamma\lambda_d\left(\bigcup_{i=1}^kB(x_i,2u'_{\gamma}) \right)} \times \mu_k(x_1,\ldots, x_K).\end{multline}
The term $\mu_k(x_1,\ldots, x_K)$ denotes the probability that $S(x_i,u'_{\gamma})$ is covered with the spherical caps $\{\mathscr{A}_{x_j}(x_i), i\neq j\leq K\}$ and $\{\mathscr{A}_{y_m}(x_i), m\leq k\}$ where $y_1,\ldots, y_k$ are $k$ independent points which are uniformly distributed in $\bigcup_{i=1}^KB(x_i,2u'_{\gamma})$ and such that $F_{1}(C_{X_{\gamma}\cup\{\mathbf{x}_K\}}(x_i))\geq 4$ for all $i\leq K$. This probability satisfies the following property:

\begin{Le}
\label{Calottesimplex}
Let $\mathbf{x}_K=(x_1,\ldots, x_K)\in W_K(n_1,\ldots, n_K)\subset\RR^2$ and \begin{equation}\label{defN'}N'=4n_1+4n_2+\sum_{l=3}^K3n_l.\end{equation} Then, for all $k<N'$ 
 \begin{equation}\label{coverKspheressimplex}\mu_k(x_1,\ldots, x_K)=0.\end{equation}  
\end{Le}
The proof of Lemma \ref{Calottesimplex} is postponed to the appendix. From \eqref{jointdistributionsimplex}, \eqref{coverKspheressimplex} and \eqref{defN'}, we deduce for $\gamma$ large enough that\[\gamma^K\PPP{\forall i\leq K, R(C_{X_{\gamma}\cup\{\mathbf{x}_K\}}(x_i))<u'_{\gamma}, F_{1}(C_{X_{\gamma}\cup\{\mathbf{x}_K\}}(x_i))\geq 4}\leq c\cdot \gamma^K(\gamma u'^2_{\gamma})^{N'} = c\cdot\gamma^{n_2}\prod_{l=3}^K\gamma^{\frac{4l-3}{4}n_l}.\]
 Moreover, $\lambda_{2K}(W_{K}(n_1,\ldots, n_K))\leq c\cdot\prod_{l=2}^K(u'^2_{\gamma})^{(l-1)n_l} = c\cdot \gamma^{-\frac{5}{4}n_2}\prod_{l=3}^K\gamma^{\frac{-5l+5}{4}n_l}$. This shows that \begin{equation}
\label{minsimplexstep2.2}
\gamma^K\int_{W^K-W_K(K,0,\ldots, 0)}\PPP{\forall i\leq K, R(C_{X_{\gamma}\cup\{\mathbf{x}_K\}}(x_i))<u'_{\gamma}, F_{1}(C_{X_{\gamma}\cup\{\mathbf{x}_K\}}(x_i))\geq 4}d\mathbf{x}_K = O\left(\gamma^{-1/4} \right)
.\end{equation} From \eqref{minsimplexstep2.1}, \eqref{minsimplexstep2.2} and Lemma \ref{rewriteHenze}, we get \[\PPP{R'_{\min}(\gamma)\geq u'_{\gamma}} \conv[\gamma]{\infty}e^{-t^{4}}.\]
\end{prooft}

We conclude the section with a quick sketch of proof for \eqref{maxins}. 

\begin{prooft}{\eqref{maxins}}
We notice that \[r_{\max}(\gamma) = \max_{x\in X_{\gamma}\cap W}r(C_{X_{\gamma}}(x)) = \frac{1}{2}\max_{x\in X_{\gamma}\cap W}\min_{y\neq x\in X_{\gamma}}d(x,y).\] The behaviour of the maximum of nearest neighbor distances was studied by Henze in Theorem 1 of \cite{He} when the input is a binomial process. His result did not include the contribution of boundary effects and is consequently limited to the set of points in $W\ominus B(0,u_{\gamma})$. With Lemma \ref{rewriteHenze} and proceeding along the same lines as in the proof of \eqref{mincircons}, we are able to show the convergence in distribution of the maximal inradius of Voronoi tessellation when the input is a Poisson point process in $W$. 
\end{prooft}

\section{Proof of \eqref{maxcircons}, consequence on Poisson-Voronoi approximation}
\label{maxcirconsparagraph}
\begin{prooft}{\eqref{maxcircons}}
First, we notice that \[R_{\max}(\gamma) = \max_{x\in X_{\gamma}\cap W}R(C_{X_{\gamma}}(x)) = \max_{x\in X_{\gamma}\cap W} \max_{y\in C_{X_{\gamma}}(x)}d(x,y).\] In order to avoid boundary effects, we start by studying an intermediary radius $R'_{\max}(\gamma)$ defined as 
\[R'_{\max}(\gamma)=\max_{x\in X_{\gamma}, C_{X_{\gamma}}(x)\cap W\neq \varnothing} \max_{y\in C_{X_{\gamma}}(x)\cap W}d(x,y).\] In a first step, we provide the asymptotic behaviour of $R'_{\max}(\gamma)$. Secondly, we study the effects of Voronoi cells astride $W$ and $W^c$.  

\paragraph{Step 1} The distribution function of $R'_{\max}(\gamma)$ can be interpreted as a covering probability. Indeed, if we denote by \begin{equation}\label{deft}u_{\gamma}=u_{\gamma}(t)=\left(\frac{1}{\kappa_d\gamma}t+\frac{1}{\kappa_d\gamma}\log\left(\alpha_1\gamma (\log\gamma)^{d-1}\right)\right)^{1/d}\end{equation} where \begin{equation}\label{defalpha1}
\alpha_1:=\frac{1}{d!}\left(\frac{\pi^{1/2}\Gamma\left(\frac{d}{2}+1\right)}{\Gamma\left(\frac{d+1}{2}\right)}\right)^{d-1}
\end{equation} and $t$ is a fixed parameter, we have \[\begin{split} R'_{\max}(\gamma) \leq u_{\gamma} &   \Longleftrightarrow \forall x\in X_{\gamma}, \text{ s.t. } C_{X_{\gamma}}(x)\cap W\neq \varnothing, \forall y\in C_{X_{\gamma}}(x)\cap W, d(x,y)\leq u_{\gamma}\\
& \Longleftrightarrow \forall y\in W, \exists x\in X_{\gamma}, d(x,y)\leq u_{\gamma}\\
& \Longleftrightarrow\left\{ B(x,u_{\gamma}),x\in X_{\gamma}\right\} \text{ covers } W. \end{split}\]
We have to deal with the probability to cover a region with a large number of balls having a small radius when $\gamma\rightarrow\infty$. Asymptotics of such covering probabilities have been studied by Janson. We apply Lemma 7.3 of \cite{J} which is rewritten in our particular framework. Actually, Lemma 7.3 of \cite{J} investigates covering with copies of a general convex body and requires conditions which are clearly satisfied in the case of the ball (see Lemmas 5.2, 5.4 and (9.24) therein). 
\begin{Le}(Janson)
\label{LemmaJanson}
Let $W$ be a bounded subset of $\RR^d$  such that $\lambda_d(\partial W)=0$ and $X_{\gamma}$ a Poisson point process of intensity $\gamma$. Let $R$ a random variable such that $\EE[R]>0$ and $\EE[R^{d+\epsilon}]$ for some $\epsilon >0$. We denote by $\alpha(B(0,R)) = \alpha_1\EE[R^{d-1}]^d\EE[R^d]^{-(d-1)}$. If $a=a(\gamma)$ is a function such  that $a(\gamma)\conv[\gamma]{\infty}0$ and \begin{equation}\label{CondJanson}\EEE{\lambda_d(aB(0,R))}\gamma - \log\frac{\lambda_d(W)}{\EEE{\lambda_d(aB(0,R))}} -d\log\log\frac{\lambda_d(W)}{\EEE{\lambda_d(aB(0,R))}} - \log\alpha(B(0,R))  \conv[\gamma]{\infty}t, -\infty < t < \infty\end{equation} then \[\PPP{\left\{B(x,R),x\in X_{\gamma}\right\} \text{ covers } W}\conv[\gamma]{\infty}e^{-e^{-t}}.\]

\end{Le}
Taking $a=u_{\gamma}$, $R=1$, $\lambda_d(W)=1$ and noting that $\EEE{\lambda_d(aB(0,R))}=\kappa_du_{\gamma}^d$ and $\alpha(B(0,R))=\alpha_1$, we check easily \eqref{CondJanson}. From Lemma \ref{LemmaJanson}, we deduce that $\PP\left(\left\{ B(x,u_{\gamma}),x\in X_{\gamma}\right\} \text{ covers } W \right)$ converges to $e^{-e^{-t}}$. Hence, for all $t\in\RR$, 
 \begin{equation}\label{maxcircons3}\lim_{\gamma\rightarrow\infty}\PP\left(R'_{\max}(\gamma)\leq u_{\gamma}\right)= e^{-e^{-t}}.\end{equation}

\paragraph{Step 2}
Taking $f(C_{X_{\gamma}}(x)) =  \kappa_d(\max_{y\in C_{X_{\gamma}}(x)\cap W}d(x,y))^d$, $a_{\gamma} = \gamma$, $b_{\gamma} = \gamma\log\left(\alpha_1\gamma (\log\gamma)^{d-1}\right)$ and $Y$ a Gumbel distribution (i.e. $\PP(Y\leq t) = e^{-e^{-t}}, t\in\RR$), one can check condition \eqref{hypothesisab} with $k=d$. From \eqref{maxcircons3} and Proposition \ref{boundaryLemma}, we deduce that $\PP(\max_{x\in X_{\gamma}\cap W}\max_{y\in C_{X_{\gamma}}(x)\cap W}d(x,y)\leq u_{\gamma})$ converges to $e^{-e^{-t}}$ for all $t\in\RR$. Using the fact that, on the event $A_{\gamma}$ (given in Lemma \ref{HSS}), \[\max_{x\in X_{\gamma}\cap W}\max_{y\in C_{X_{\gamma}}(x)\cap W}d(x,y) \leq \max_{x\in X_{\gamma}\cap W}\max_{y\in C_{X_{\gamma}}(x)}d(x,y) \leq \max_{x\in X_{\gamma}\cap W_{1+l_{\gamma}}}\max_{y\in C_{X_{\gamma}}(x)\cap W_{1+l_{\gamma}}}d(x,y)\]
 and proceeding along the same lines as in the proof of Proposition \ref{boundaryLemma}, we get \begin{equation}\label{maxcircons5}\PPP{R_{\max}(\gamma)\leq u_{\gamma}} \conv[\gamma]{\infty}e^{-e^{-t}}.\end{equation}
\end{prooft}
We can note that the asymptotic behaviour of $R_{\max}(\gamma)$ gives an interpretation of Lemma 7.3  in \cite{J}. Indeed, \eqref{maxcircons5} shows that the Gumbel distribution which appears as a limit probability of a covering is actually the limit distribution of a maximum.

We now apply this convergence result to the so-called Poisson-Voronoi approximation defined as 
\[\mathscr{V}_{X_{\gamma}}(W)=\bigcup_{x\in X_{\gamma}}C_{X_{\gamma}}(x).\] It consists in discretizing a given convex window $W$ with a finite union of convex polyhedra. This approximation has various applications such as image analysis (reconstructing an image from its intersection with a Poisson point process, see \cite{KT}) or quantization (see chapter 9 of \cite{GL}). Estimates of the first two moments of the symmetric difference between the convex body and its approximation are given in \cite{HR} and extended to higher moments in \cite{RSZ}. To the best of our knowledge, the convergence of $\mathscr{V}_{X_{\gamma}}(W)$ to $W$ in the sense of Hausdorff distance, denoted by $d_H(\cdot, \cdot)$,  has not been investigated. Corollary \ref{Hausdorff} addresses that question with an assumption on the regularity of $W$ which is in the spirit of the $n$-regularity (see Definition 3 in \cite{CV}).

\begin{Cor}
\label{Hausdorff}
Let us assume that there exists $\alpha >0$ such that, for $v$ small enough and for all $y\in W$, \begin{equation}\label{hyphausdorff}\lambda_d(B(y,v)\cap W)\geq \alpha\lambda_d(B(y,v)).\end{equation} Then 
  \begin{equation}\label{Hausdorffspeed}\PPP{d_H(W,\mathscr{V}_{X_{\gamma}}(W))\leq \left(c(\alpha)\gamma^{-1}\log\left(\alpha_1\gamma(\log\gamma)^{d-1} \right)\right)^{1/d}}\conv[\gamma]{\infty}1\end{equation}  where $c(\alpha)=\kappa_d^{-1}+2^d\kappa_d^{-1}\alpha^{-1}$. 
\end{Cor}

\begin{prooft}{Corollary \ref{Hausdorff}}
Let us denote by \begin{equation}\label{defvgamma}v_{\gamma} = \left(c(\alpha){\gamma}^{-1}\log\left(\alpha_1\gamma(\log\gamma)^{d-1} \right)\right)^{1/d}.\end{equation} 
First, we show that $\max_{y\in\mathscr{V}_{X_{\gamma}}(W)}d(y,W)\leq v_{\gamma}$ with high probability. For all $t\in\RR$, using the fact that $u_{\gamma}\leq v_{\gamma}$ for $\gamma$ large enough, where $u_{\gamma}=u_{\gamma}(t)$ is given in \eqref{deft}, we get \[\PPP{\max_{y\in\mathscr{V}_{X_{\gamma}}(W)}d(y,W)\leq v_{\gamma}} \geq \PPP{R_{\max}(\gamma)\leq v_{\gamma}} \geq \PPP{R_{\max}(\gamma)\leq u_{\gamma}}.\] From \eqref{maxcircons5} and  Proposition \ref{boundaryLemma}, the last term converges to $e^{-e^{-t}}$ as $\gamma$ goes to infinity. Taking $t\rightarrow\infty$, we get \begin{equation}\label{Hausdorffpart1}
\lim_{\gamma\rightarrow\infty}\PPP{\max_{y\in\mathscr{V}_{X_{\gamma}}(W)}d(y,W)\leq v_{\gamma}}\geq \lim_{t\rightarrow\infty}e^{-e^{-t}}=1
.\end{equation}

In a second step, we are going to show that $\max_{y\in W}d(y,X_{\gamma}\cap W)\leq v_{\gamma}$ with high probability via the use of a covering of $W$ by balls as in the proof of \eqref{maxcircons}. Now, the convex body $W$ is covered by $\mathscr{N} = O\left(v_{\gamma}^{-d} \right)$ deterministic balls $B_1,\ldots, B_{\mathscr{N}}$ with center in $W$ and radius equal to $v_{\gamma}/2$. From \eqref{hyphausdorff}, \eqref{defvgamma} and the fact that $\#(B_i\cap (X_{\gamma}\cap W))$ is Poisson distributed with mean  $\gamma\lambda_d(B_i\cap W)$, we get for $\gamma$ large enough
\[\PPP{\max_{y\in W}d(y,X_{\gamma}\cap W)> v_{\gamma}} \leq \PPP{\bigcup_{i=1}^{\mathscr{N}} \{\#(B_i\cap (X_{\gamma}\cap W))=0\}} \leq \mathscr{N}e^{-\gamma\alpha\kappa_d \left(v_{\gamma}/2 \right)^d}\leq \alpha_1^{-1}\gamma^{-1}(\log\gamma)^{-(d-1)}\mathscr{N}.\] Using the fact that $\mathscr{N}=O(v_{\gamma}^{-d})$ i.e. $\mathscr{N}=O(\gamma(\log\gamma)^{-1})$ according to \eqref{defvgamma}, the right-hand side is $O\left((\log\gamma)^{-d}\right)$. Hence 
\begin{equation}\label{Hausdorffpart2}
\PPP{\max_{y\in W}d(y,X_{\gamma}\cap W)\leq v_{\gamma}}\conv[\gamma]{\infty}1
.\end{equation}
Since $d_H(W,\mathscr{V}_{X_{\gamma}}(W)) \leq \max\left\{\max_{y\in\mathscr{V}_{X_{\gamma}}(W)}d(y,W), \max_{y\in W}d(y,X_{\gamma}\cap W)\right\}$, we deduce from  \eqref{Hausdorffpart1} and \eqref{Hausdorffpart2} that \[\PPP{d_H(W,\mathscr{V}_{X_{\gamma}}(W))\leq v_{\gamma}}\conv[\gamma]{\infty}1.\]
\end{prooft}

In \cite{HR}, Heveling and Reitzner obtain that the volume of the symmetric difference between $W$ and $\mathscr{V}_{X_{\gamma}}(W)$ is of the order of $\gamma^{-1/d}$. The result above makes sense and could provide the right order of the Hausdorff distance. Obviously, the constant $c(\alpha)=\kappa_d^{-1}+2^d\kappa_d^{-1}\alpha^{-1}$ is not optimal. From Lemma \ref{HSS}, it would have been possible to get an upper-bound of the order of $\gamma^{-(1-\epsilon)/d}$ but it is less precise than Corollary \ref{Hausdorff}. 

\section{Proof of \eqref{minins}}
\begin{prooft}{\eqref{minins}}
Let $t\geq 0$ be fixed. We denote by $u_{\gamma}$ the following function: \begin{equation}\label{defugammaminins}u_{\gamma} = u_{\gamma}(t) = \left(2^{-(d-1)}\kappa_d^{-1}\gamma^{-2}t \right)^{1/d}.\end{equation} We start by finding a different expression of $r_{\min}(\gamma)$ which does not rely on the Voronoi structure. Indeed, for all $x\in X_{\gamma}\cap W$ we have \[r(C_{X_{\gamma}}(x)) = \max\{r\geq 0, B(x,r)\subset C_{X_{\gamma}}(x)\} = \frac{1}{2}\min_{y\neq x\in X_{\gamma}}d(x,y).\] Hence, $r_{\min}(\gamma)$ can be rewritten as \begin{equation}\label{rewriteminins}r_{\min}(\gamma) = \frac{1}{2}\min_{(x,y)_{\neq}\in (X_{\gamma}\cap W)\times X_{\gamma}}d(x,y).
\end{equation}  The equality \eqref{rewriteminins} implies that the problem is reduced to a study of inter-point distance. Such study is well known for a binomial process $X^{(n)}$ of intensity $n$ in $W$. In particular, Jammalamadaka and Janson (see \cite{JJ}, §4) have shown that for all $t\geq 0$,  
\begin{equation}
\label{convbinom}
\PPP{r'_{\min, n}\geq u_n}\conv[n]{\infty}e^{-t}
\end{equation} where $r'_{\min,n}$ is defined as \[r'_{\min,n} = \frac{1}{2}\min_{(x,y)_{\neq}\in X^{(n)}\times X^{(n)}}d(x,y)\] and $u_n$ given in \eqref{defugammaminins}.  In a first elementary step, we extend the limit to a Poisson point process. Our main contribution is then to compare the obtained limit with $r_{\min}(\gamma)$ by dealing with boundary effects. In particular, our study provides a far more accurate estimate of the contribution of boundary cells (see \eqref{convexp}) than what we could have deduced from Proposition \ref{boundaryLemma}.  

\paragraph{Step 1} We extend \eqref{convbinom} to a Poisson point process. We define  \begin{equation}\label{defins}r'_{\min}(\gamma) = \frac{1}{2}\min_{(x,y)_{\neq}\in (X_{\gamma}\cap W)^2}d(x,y).\end{equation}
 Let us note that for all $0<\alpha<\beta<1$, and for all $n\in\{0,1,2,\ldots\}$, $|n-\gamma|\leq \gamma^{\alpha}\Longrightarrow |n-\gamma|\leq n^{\beta}$ for $\gamma$ large enough. Consequently, since $u_{\gamma}$ is non-increasing in $\gamma$, we have for $\gamma$ large enough

\begin{multline}
\label{insbin}
\left|\PP\left(r'_{\min}(\gamma)\geq u_{\gamma} \right) - e^{-t} \right|  \leq \sum_{n=0}^{\infty}\left|\PP(r'_{\min,n}\geq u_{\gamma}) - e^{-t} \right|\PP(\#(X_{\gamma}\cap W)=n)\\
 \leq \sum_{|n-\gamma|\leq \gamma^{\alpha}}\max\left\{\left|\PP(r'_{\min,n}\geq u_{n-n^{\beta}}) - e^{-t} \right|, \left|\PP(r'_{\min,n}\geq u_{n+n^{\beta}}) - e^{-t} \right|\right\}\PP(\#(X_{\gamma}\cap W)=n)\\ + \PP(|\#(X_{\gamma}\cap W)-\gamma|>\gamma^{\alpha})
.\end{multline} The second term of \eqref{insbin} converges to 0 thanks to a concentration inequality for Poisson variables (see e.g. Lemma 1.4 in \cite{Pr}). The first term is lower than $\max_{n\geq \gamma - \gamma^{\alpha}}\max\left\{\left|\PP(r'_{\min,n}\geq u_{n-n^{\beta}}(t)) - e^{-t} \right|, \left|\PP(r'_{\min,n}\geq u_{n+n^{\beta}}(t)) - e^{-t} \right|\right\}$ which tends to 0 according to \eqref{convbinom}. This shows that, for all $t\geq 0$, 
 \begin{equation}\label{inspoiss}\lim_{\gamma\rightarrow\infty}\PP\left(r'_{\min}(\gamma)\geq u_{\gamma} \right) = e^{-t}.\end{equation}

\paragraph{Step 2} We show that $r_{\min}(\gamma)=r'_{\min}(\gamma)$ with probability of order of $O(\gamma^{-\epsilon})$ with $\epsilon\in (0,\frac{2}{d})$. Indeed, the random variables $r_{\min}(\gamma)$ and $r'_{\min}(\gamma)$, defined in \eqref{rewriteminins}  and \eqref{defins}, are equal if and only if no point of $X_{\gamma}\cap W^c$ falls into the union of the balls $B(x,2r'_{\min}(\gamma))$ for $x\in X_{\gamma}\cap W$ such that $d(x,\partial W)<2r'_{\min}(\gamma)$ i.e. 
\begin{multline}
\label{majexp}
\PP(r_{\min}(\gamma) \neq r'_{\min}(\gamma)) = \PP\left(\#\left( X_{\gamma}\cap W^c\cap \bigcup_{\underset{d(x,\partial W)<2r'_{\min}(\gamma)}{x\in X_{\gamma}\cap W,}}B(x,2r'_{\min}(\gamma))\right) \neq 0\right)\\\leq  \EE\left[\sum_{\underset{d(x,\partial W)<2r'_{\min}(\gamma)}{x\in X_{\gamma}\cap W,}}\#\left(X_{\gamma}\cap W^c\cap B(x,2r'_{\min}(\gamma))\right)\right]
.\end{multline} From Slivnyak-Mecke formula (see e.g. Corollary 3.2.3 of \cite{SW}), we get
\[\EE\left[\sum_{\underset{d(x,\partial W)<2r'_{\min}(\gamma)}{x\in X_{\gamma}\cap W,}}\#\left(X_{\gamma}\cap W^c\cap B(x,2r'_{\min}(\gamma))\right)\right] = \int_{W}\gamma \EE\left[\#\left(X_{\gamma}\cap W^c\cap B(x,2r'^{(x)}_{\min}(\gamma))\right) \mathbb{1}_{d(x,\partial W)<2r'^{(x)}_{\min}(\gamma)}\right]dx\] where $r'^{(x)}_{\min}(\gamma) = \frac{1}{2}\min_{(x',y)_{\neq}\in (X_{\gamma}\cup\{x\}\cap W)^2}d(x',y) $ for all $x\in X_{\gamma}\cap W$. Noting that $r'^{(x)}_{\min}(\gamma)\leq r'_{\min}(\gamma)$, we then obtain \begin{multline}\label{majexp2} \EE\left[\sum_{\underset{d(x,\partial W)<2r'_{\min}(\gamma)}{x\in X_{\gamma}\cap W,}}\#\left(X_{\gamma}\cap W^c\cap B(x,2r'_{\min}(\gamma))\right)\right] \leq \int_{W}\gamma \EE\left[\#\left(X_{\gamma}\cap W^c\cap B(x,2r'_{\min}(\gamma))\right) \mathbb{1}_{d(x,\partial W)<2r'_{\min}(\gamma)}\right]dx\\
 = \int_{W}\gamma \EEE{\EE\left[\#\left(X_{\gamma}\cap W^c\cap B(x,2r'_{\min}(\gamma))\right) |X_{\gamma}\cap W\right]\mathbb{1}_{d(x,\partial W)<2r'_{\min}(\gamma)}}dx   .\end{multline} Since  $\#\left(X_{\gamma}\cap W^c\cap B(x,2r'_{\min}(\gamma))\right)$ is Poisson distributed, we get  
\begin{multline}\label{majexp3}\gamma\EEE{\EE\left[\#\left(X_{\gamma}\cap W^c\cap B(x,2r'_{\min}(\gamma))\right) |X_{\gamma}\cap W\right]\mathbb{1}_{d(x,\partial W)<2r'_{\min}(\gamma)}}\\ 
 = \gamma^2\EE\left[\lambda_d(W^c\cap B(x,2r'_{\min}(\gamma)))\mathbb{1}_{d(x,\partial W)<2r'_{\min}(\gamma)} \right] \leq  2^d\kappa_d\cdot  \gamma^2\EE\left[{r'_{\min}(\gamma)}^d\mathbb{1}_{d(x,\partial W)<2r'_{\min}(\gamma)}\right] . \end{multline} Using \eqref{majexp}, \eqref{majexp2}, \eqref{majexp3} and  Fubini's theorem, we obtain 
 \begin{equation}\label{boundeffectineq1}\PP(r_{\min}(\gamma) \neq r'_{\min}(\gamma))\leq 2^d\kappa_d\cdot \gamma^2\EEE{ r'^d_{\min}(\gamma)\int_{W}\mathbb{1}_{d(x,\partial W)<2r'_{\min}(\gamma)}dx}\leq c\cdot \gamma^2\EEE{r'^{d+1}_{\min}(\gamma)}.\end{equation} The last inequality comes from Steiner formula (see (14.5) in \cite{SW}) and $c$ denotes a constant depending on $W$. Hence, to show that \begin{equation}\label{aimstep2minins}\PP\left(r_{\min}(\gamma)\neq r'_{\min}(\gamma) \right)\conv[\gamma]{\infty}0\end{equation} we have to find some upper-bound of $\gamma^2\EE[r'^{d+1}_{\min}(\gamma)]$. We know, from \eqref{inspoiss} and \eqref{defugammaminins} , that $\gamma^2r'^{d}_{\min}(\gamma)$ tends to 0 in distribution but it does not imply \eqref{aimstep2minins}. Lemma \ref{discr} below provides an estimate of the deviation probabilities of $\gamma^2r'^d_{\min}(\gamma)$ when the window W is a cube.

\begin{Le}
\label{discr}
Let $C$ be a cube of side $M$ and $X_{\gamma}$ a Poisson point process of intensity $\gamma$. Let us denote by \[r'_{\min\mid C}(\gamma)=\frac{1}{2}\min_{(x,y)_{\neq}\in (X_{\gamma}\cap C)^2}d(x,y).\] Then, for all $u\leq \min\{\frac{1}{4}Md^{1/2}, \frac{1}{2}d^{1/2}\gamma^{-1/d}\}$, there exists a constant $c(M)$ such that
\[\PP(r'_{\min\mid C}(\gamma)\geq u ) \leq e^{-c(M)\gamma^2u^d}.\] 
\end{Le}

\begin{prooft}{Lemma \ref{discr}}
Let $u\leq \min\{\frac{1}{4}Md^{1/2}, \frac{1}{2}d^{1/2}\gamma^{-1/d}\}$ be fixed.

We subdivide the cube $C=[0,M]^d$ into a set of $\mathscr{N}$ subcubes $C_1, \ldots, C_{\mathscr{N}}$ of equal size $c$ with $c = 2d^{-1/2}u$ and $\mathscr{N} =\left\lfloor Mc^{-1}\right\rfloor^d $. Since $\diam(C_i) = 2u$ for each $i\leq \mathscr{N}$, we obtain 
\[\PP(r'_{\min\mid C}(\gamma)\geq u) \leq \PP\left(\bigcap_{i=1}^{\mathscr{N}}\{\# (C_i\cap X_{\gamma}) \leq 1\} \right) = \left(e^{-\gamma c^d}(1+\gamma c^d)\right)^{\mathscr{N}} .\] Replacing $c^d$ by $2^dd^{-d/2}u^d$ and $\mathscr{N}$ by $\lfloor 2^{-1}Md^{1/2}u^{-1}\rfloor^d$ we obtain the following inequality:
 \[\begin{split}\PP(r'_{\min\mid C}(\gamma)\geq u) & \leq e^{{\lfloor 2^{-1}Md^{1/2}u^{-1}\rfloor}^d\left(\log(1+\gamma 2^dd^{-d/2}u^d)-\gamma 2^dd^{-d/2}u^d \right)}\\
 & \end{split}.\] 
Since $\gamma 2^dd^{-d/2}u^d\leq 1$ and $2M^{-1}d^{-1/2}u\leq \frac{1}{2}$, we have $\log(1+\gamma 2^dd^{-d/2}u^d)-\gamma 2^dd^{-d/2}u^d \leq -\frac{1}{4}2^{2d}d^{-d}\gamma^2u^{2d}$ and ${\lfloor 2^{-1}Md^{1/2}u^{-1}\rfloor}^d \geq (2^{-1}Md^{1/2}u^{-1}-1)^d\geq 2^{-2d}M^dd^{d/2}u^{-d}$. Hence \[\PP(r'_{\min\mid C}(\gamma)\geq u) \leq e^{-\frac{1}{4} d^{-d/2}M^d\gamma^2u^d} = e^{-c(M)\gamma^2u^d}\] where $c(M) = \frac{1}{4} d^{-d/2}M^d$.
\end{prooft}

Now, we can derive an upper-bound of $\gamma^2\EE[r'^{d+1}_{\min}(\gamma)]$. Indeed, since $W$ has non-empty interior, there exists a cube $C$ of side $M$ included in $W$. Using the fact that  $\#(X_{\gamma}\cap C)\geq 2 \Longrightarrow r'_{\min}(\gamma)\leq r'_{\min\mid C}(\gamma)$, we get \begin{multline}\label{boundeffecteq2}\gamma^2\EE[r'^{d+1}_{\min}(\gamma)] = \gamma^2\int_0^{\diam(W)}\PP(r'^{d+1}_{\min}(\gamma)\geq s)ds \\\leq \diam(W)\gamma^2\PP(\#(X_{\gamma}\cap C)\leq 1)+\gamma^2\int_{0}^{Md^{1/2}}\PP(r'^{d+1}_{\min\mid C}(\gamma)\geq s)ds .\end{multline} The first term of the right-hand side of \eqref{boundeffecteq2} is decreasing exponentially fast to 0 since $\#(X_{\gamma}\cap C)$ is Poisson distributed of mean $\gamma M^d$. For the second term,  let us consider a fixed $\epsilon$ in $(0,\frac{2}{d})$. Then \begin{equation}\begin{split}\label{majmin2} \gamma^2\int_{0}^{Md^{1/2}}\PP(r'^{d+1}_{\min\mid C}(\gamma)\geq s)ds & = \int_0^{\gamma^{-(2+\epsilon)}}\gamma^2\PP\left(r'_{\min\mid C}(\gamma)\geq s^{1/(d+1)}\right)ds + \int_{\gamma^{-(2+\epsilon)}}^{Md^{1/2}}\gamma^2\PP\left(r'_{\min\mid C}(\gamma)\geq s^{1/(d+1)}\right)ds\\
 & \leq \gamma^{-\epsilon}+Md^{1/2} \gamma^2\PP\left(r'_{\min\mid C}(\gamma)\geq \gamma^{-(2+\epsilon)/(d+1)}\right).
 \end{split}\end{equation} Since $\epsilon>0$,  we have $\gamma^{-(2+\epsilon)/(d+1)}\leq \min\{\frac{1}{4}Md^{1/2}, \frac{1}{2}d^{1/2}\gamma^{-1/d}\}$ for $\gamma$ large enough. Hence, from Lemma \ref{discr} applied to $u:=\gamma^{-(2+\epsilon)/(d+1)}$, we deduce that for $\gamma$ large enough, 
 \begin{equation}\label{majintmin}\gamma^2\int_{0}^{Md^{1/2}}\PP(r'^{d+1}_{\min\mid C}(\gamma)\geq s)ds\leq \gamma^{-\epsilon} + Md^{1/2}\gamma^2e^{-c(M)\gamma^{(2-\epsilon d)/(d+1)}}.\end{equation}  The last term of the right-hand side of \eqref{majintmin} converges exponentially fast to 0 as $\gamma$ goes to infinity since $\epsilon < \frac{2}{d}$. Combining that argument with \eqref{boundeffectineq1}, \eqref{boundeffecteq2} and \eqref{majintmin}, we deduce that  \begin{equation}\label{convexp}\PPP{r_{\min}(\gamma) \neq r'_{\min}(\gamma)} = O\left(\gamma^{-\epsilon} \right).\end{equation} We then deduce from \eqref{inspoiss} and \eqref{convexp} that \begin{equation*}\left|\PP\left(r'_{\min}(\gamma) \geq u_{\gamma} \right) - e^{-t} \right| \leq \left|\PP\left(r'_{\min}(\gamma) \geq u_{\gamma}  \right) - e^{-t} \right| + 2 \PP\left(r_{\min}(\gamma)\neq r'_{\min}(\gamma) \right) \underset{\gamma\rightarrow\infty}{\longrightarrow}0.\end{equation*} 
 \end{prooft}

\begin{Rk}
The rate for the convergence in distribution of $r_{\min}(\gamma)$ to the Weibull distribution can be estimated. For instance, we can show that Theorem 2.1 in \cite{ST} implies the rate of convergence of $r'_{\min}(\gamma)$. Another way to get it is to use  Theorem 1 in \cite{AGG}. We then obtain that there exists positive constants $c(W)$ and $\Gamma(W)$ such that, for all $\epsilon<\frac{2}{d}$, $t\geq 0$ and $\gamma\geq \Gamma(W)$, \[\left|\PP\left(2^{d-1}\kappa_d\gamma^2 r_{\min}(\gamma)^d \geq t\right) - e^{-t}\right| \leq c(W)\gamma^{-\min\{\frac{1}{2},\epsilon\}}.\]
\end{Rk}

\vspace{0.5cm}

The study of more extremes for general tessellations and their rates of convergence will be developed in a future paper.

\section*{Appendix}

\begin{prooft}{Lemma \ref{Calotte}}
Actually, we show the following deterministic result: let $K\geq 2$, $k<N$, $(x_1,\ldots, x_K)\in W_{K}(n_1,\ldots, n_K)$ with $(n_1,\ldots, n_K)\neq (K,0,\ldots, 0)$ and $(y_1,\ldots, y_k)\in \bigcup_{i=1}^KB(x_i,2u_{\gamma})$ such that $\{\mathbf{x}_K\}\cup\{\mathbf{y}_k\}$ are in general position i.e. each subset of size $n<d+1$ is affinely independent (see \cite{Ze}). Then there exists $i\leq K$ such that sphere $S(x_i,u_{\gamma})$ is not covered by the induced spherical caps $\{\mathscr{A}_{x_j}(x_i), i\neq j\leq K\}\cup\{\mathscr{A}_{y_m}(x_i), m\leq k\}$. 

Indeed, from \eqref{defN} there exists a connected component of $\bigcup_{i=1}^KB(x_i,2u_{\gamma})$ of size $1\leq l\leq K$, say $\bigcup_{i=1}^lB(x_{i},2u_{\gamma})$ without loss of generality, such that $N_l< d+1$ with
\begin{equation}\label{defNl}N_l = \#\left(\{\mathbf{y}_k\}\cap\bigcup_{i=1}^lB(x_i,2u_{\gamma}) \right)\end{equation}
 
Since $\{\mathbf{x}_l\}\cup\{\mathbf{y}_{N_l}\}$ are in general position, the family $\{\mathbf{x}_l\}$ is not included in the convex hull of $\{\mathbf{y}_{N_l}\}$. In particular, there exists $i\leq l$ such that $x_i$ is not in the convex hull of $\{\mathbf{x}_l\}\cup\{\mathbf{y}_{N_l}\}-\{x_i\}$. Since a Voronoi cell induced by a finite number of points is not bounded if and only if its nucleus is an extremal point of the polytope induced by the points, it implies that the circumscribed radius of  $C_{\{\mathbf{x_l}\}\cup\{\mathbf{y}_{N_l}\}}(x_i)$ is not finite i.e. $S(x_i,u_{\gamma})$ is not covered. 
\end{prooft}

\begin{prooft}{Lemma \ref{Calottesimplex}}
We show the following deterministic result: let $K\geq 2$, $k<N'$, $(x_1,\ldots, x_K)\in W_K(n_1,\ldots, n_K)$ with $(n_1,\ldots, n_K)\neq (K,0,\ldots, 0)$ and $(y_1,\ldots, y_k)\in \bigcup_{i=1}^KB(x_i,2u_{\gamma})$ such that $\{\mathbf{x}_K\}\cup\{\mathbf{y}_k\}$ are in general position. Then there exists $i\leq K$ such that either the sphere $S(x_i,u_{\gamma})$ is not covered by the induced spherical caps $\{\mathscr{A}_{x_j}(x_i), i\neq j\leq K\}\cup\{\mathscr{A}_{y_m}(x_i), m\leq k\}$ or $F_1(C_{\{\mathbf{x}_K\}\cup\{\mathbf{y}_k\}}(x_i)) \leq 3$. 

Indeed, from \eqref{defN'} there exists a connected component of $\bigcup_{i=1}^KB(x_i,2u_{\gamma})$ of size $1\leq l\leq K$, say $\bigcup_{i=1}^lB(x_{i},2u_{\gamma})$ without loss of generality, such that $N_l< 4$ if $l=1,2$ and $N_l<3$ if $l\geq 3$ where $N_l$ is given in \eqref{defNl}. 
\begin{itemize}
\item If $l=1$, either $S(x_1,u_{\gamma})$ is covered or $F_1(C_{\{\mathbf{x}_K\}\cup\{\mathbf{y}_k\}}(x_1)) = F_1(C_{\{x_1\}\cup\{\mathbf{y}_{N_1}\}}(x_1)) \leq 3$ since $N_1\leq 3$. 
\item If $l\geq 3$, from Lemma \ref{Calotte} there exists $i\leq l$ such that $S(x_i,u_{\gamma})$ is not covered.
\item If $l=2$, we can assume that $N_2=3$. We have to prove that if $\mathbf{y}_3=\{y_1,y_2,y_3\}$ is a set of three points in $B(x_1,2u_{\gamma})\cup B(x_2,2u_{\gamma})$, then the following properties 1 and 2 below cannot hold simultaneously.
\begin{enumerate}
\item The circles $S(x_1,u_{\gamma})$ and $S(x_2,u_{\gamma})$ are covered by the induced circular caps \[\{\mathscr{A}_{x_1}(x_2), \mathscr{A}_{x_2}(x_1)\}\cup\{\mathscr{A}_{y_m}(x_i), m\leq 3\}.\] 
\item The number of edges of the Voronoi cells satisfy $F_1(C_{\{x_1,x_2,\mathbf{y}_3\}}(x_1))\geq 4$ and $F_1(C_{\{x_1,x_2,\mathbf{y}_3\}}(x_2))\geq 4$. 
\end{enumerate}
Let us assume that Properties 1 and 2 hold simultaneously. Let us denote by $G$ the Delaunay graph associated to $\{x_1,x_2,y_1,y_2,y_3\}$. Then $G$ is a connected planar graph with $v=5$ vertices and $e$ edges. From Euler's formula on planar graphs, $e\leq 3v-6$ i.e. \begin{equation}\label{condedge}
e\leq 9
.\end{equation} From Property 1 and according to the proof of Lemma \ref{Calotte}, $x_1,x_2$ are in the convex hull of $\{y_1,y_2,y_3\}$ i.e. $\{y_1,y_2\}$, $\{y_1,y_3\}$ and $\{y_2,y_3\}$ are edges of the associated Delaunay triangulation. From Property 2, $x_1,x_2$ are connected to every point i.e. $\{x_1,x_2\}$, $\{x_1,y_1\}$, $\{x_1,y_2\}$, $\{x_1,y_3\}$, $\{x_2,y_1\}$, $\{x_2,y_2\}$ and $\{x_2,y_3\}$ are also edges  of the Delaunay triangulation. The total number of these edges is $e=10$. This contradicts \eqref{condedge}. 

\end{itemize}

  \end{prooft}

\begin{center}
\textbf{Acknowledgement }
\end{center}
This work was partially supported by the French ANR grant PRESAGE (ANR-11-BS02-003) and the French research group GeoSto (CNRS-GDR3477).

\bibliographystyle{abbrv}
\bibliography{BiblioVE}

\end{document}